\RequirePackage{fix-cm}
\documentclass[twocolumn]{svjour3}          
\smartqed  
\usepackage{graphicx}
 \usepackage{mathptmx}      
\usepackage{amsmath}
\usepackage{amssymb}
\usepackage{lineno}
\usepackage{natbib}
\begin{document}

\title{Data-driven prediction strategies for low-frequency patterns of North Pacific climate variability}


\author{Darin Comeau         \and
	    Zhizhen Zhao	\and
            Dimitrios Giannakis \and
            Andrew J. Majda
}


\institute{D. Comeau \at
              Center for Atmosphere Ocean Science, Courant Institute of Mathematical Sciences, New York University, New York, NY.\\
              \email{comeau@cims.nyu.edu}           
}

%
%

\date{Received: date / Accepted: date}

\maketitle

\begin{abstract}
The North Pacific exhibits patterns of low-frequency variability on the intra-annual to decadal time scales, which manifest themselves in both model data and the observational record, and prediction of such low-frequency modes of variability is of great interest to the community. While parametric models, such as stationary and non-stationary autoregressive models, possibly including external factors, may perform well in a data-fitting setting, they may perform poorly in a prediction setting. Ensemble analog forecasting, which relies on the historical record to provide estimates of the future based on past trajectories of those states similar to the initial state of interest, provides a promising, nonparametric approach to forecasting that makes no assumptions on the underlying dynamics or its statistics. We apply such forecasting to low-frequency modes of variability for the North Pacific sea surface temperature and sea ice concentration fields extracted through Nonlinear Laplacian Spectral Analysis. We find such methods may outperform parametric methods and simple persistence with increased predictive skill.
\end{abstract}


\section{Introduction}
\label{sec:intro}

	\par Predictability in general circulation models (GCMs) for the North Pacific, from seasonal to decadal time scales, has been the subject of many recent studies, (e.g., \citet{tietsche2014seasonal}, \citet{blanchard2014characteristics}, \citet{blanchard2011persistence}), several of which focus on the role of the initial state (e.g., \citet{collins2002climate}, \citet{blanchard2011influence}, \citet{branstator2012systematic}, \citet{day2014pan}). The North Pacific exhibits prominent examples of interannual to decadal variability, such as the Pacific Decadal Oscillation \citep[PDO;][]{mantua1997pacific}, and the more rapidly decorrelating North Pacific Gyre Oscillation \citep[NPGO;][]{di2008north}, both of which have been the subject of much interest. An important phenomenon of intra-annual variability in the North Pacific is the reemergence of anomalies in both the sea surface temperature (SST) fields \citep{alexander1999reemergence}, as well as in the sea ice concentration (SIC) field \citep{blanchard2011persistence}, where regional anomalies in these state variables vanish over a season, and reappear several months later, as made evident by high time-lagged correlations. 

	\par The North Pacific (along with the North Atlantic) is a region of relative strong low-frequency variability in the global climate system \citep{branstator2012systematic}. Yet, in GCMs it has been shown that this region shows relative lack of predictability \citep[less than a decade;][]{collins2002climate}, with the Community Climate System Model (CCSM) having particularly weak persistence and low predictability in the North Pacific among similiar GCMs \citep{branstator2012systematic}. The ocean and sea ice systems show stronger low-frequency variability than the atmosphere \citep{newman2003enso}. The internal variability exhibited in the North Pacific has also been characterized as being in distinct climate regimes (e.g. \citet{overland2006regime,overland2008north}), where the dynamics exhibit regime transitions and metastability. As a result, cluster based methods have been a popular approach to model regime behavior in climate settings, such as in \citet{franzke2008hidden, franzke2009systematic}, where hidden Markov models were used to model atmospheric flows.
	
	\par A traditional modeling approach for the PDO has been to fit an autoregressive process \citep{hasselmann1976stochastic, frankignoul1977stochastic}, as well as models being externally forced from the tropics through ENSO \citep{newman2003enso}. In \citet{giannakis2012limits}, autoregressive models were successful in predicting temporal patterns corresponding to the PDO and NPGO when forced with suitably modulated intermittent modes. Additional flexibility can be built into regression models by allowing nonstationary, i.e. time dependent coefficients, and a successful example of this is the finite element method (FEM) clustering procedure combined with multivariate autoregressive factor (VARX) model framework \citep{horenko2010clustering, horenko2010identification}. In this approach, the data is partitioned into a predetermined number of clusters, and model regression coefficients are estimated on each cluster, together with a cluster affiliation function that indicates which model is used at a particular time. This can further be adapted to account for external factors \citep{horenko2011nonstationarity}. While FEM-VARX methods can be effective at fitting the desired data, advancing the system state in the future for a prediction is dependent on being able to successfully predict the unknown cluster affiliation function. Methods for using this framework in a prediction setting have been used in \citet{horenko2011nonstationarity, horenko2011analysis}. Another regression modeling approach was put forward in \citet{majda2013physics}, where physics constraints were imposed on multilevel nonlinear regression models, preventing ad-hoc, finite-time blow-up, a pathological behavior previously shown to exist in such models without physics constraints in \citet{majda2012fundamental}. Appropriate nonlinear regression models using this strategy have been shown to have high skill for predicting the intermittent cloud patterns of tropical intra-seasonal variability \citep{chen2014predicting,chen2015predicting,chen2015predictingmcw}.
	
	\par Parametric models may perform well for fitting, but often have poor performance in a prediction setting, particularly in systems that exhibit distinct dynamical regimes. Nonparametric models can be advantageous in systems where the underlying dynamical system is unknown or imperfect. An early example of this is an analog forecast, first introduced by \citet{lorenz1969atmospheric}, where one considers the historical record, and makes predictions based on examining state variable trajectories in the past that are similiar to the current state. Analog forecasting does not make any assumptions on the underlying dynamics, and cleverly avoids model error when the underlying model is observations from nature. This has since been applied to other climate prediction scenarios, such as the Southern Oscillation Index \citep{drosdowsky1994analog}, the Indian summer monsoon \citep{xavier2007analog}, and wind forecasting \citep{alessandrini2015novel}, where it was found to be particularly useful in forecasting rare events.
	
	\par Key to the success of an analog forecasting method is the ability to identify a good historical analog to the current initial state. In climate applications, the choice of analog is usually determined by minimizing Euclidean distance between snapshots of system states, with a single analog being selected \citep{lorenz1969atmospheric, branstator2012systematic}. In \citet{zhao2014analog}, analog forecasting was extended upon in two key ways. First, the state vectors considered were in Takens lagged embedding space \citep{takens1981detecting}, which captures some of the dynamics of the system, rather than a snapshot in time. Second, instead of selecting a single analog determined by Euclidean distance, weighted sums of analogs were considered, where weights are determined by a kernel function. In this context, a kernel is an exponentially decaying pairwise similarity measure, intuitively, playing the role of a local covariance matrix. In \citet{zhao2014analog}, kernels were introduced in the context of Nonlinear Spectral Analysis \citep[NLSA;][]{giannakis2012nonlinear, giannakis2013nonlinear,giannakis2014data} for decomposing high-dimensional spatiotemporal data, leading naturally to a class of low-frequency observables for prediction through kernel eigenfunctions. In \citet{bushuk2014reemergence}, the kernel used in the NLSA algorithm was adapted to be multivariate, allowing for multiple variables, possibly of different physical units, to be considered jointly in the analysis.
	
	\par The aim of this study is to develop low dimensional, data-driven models for prediction of dominant low-frequency climate variability patterns in the North Pacific, combining the approaches laid out in \citet{zhao2014analog} and \citet{bushuk2014reemergence}. We invoke a prediction approach that is purely statistical, making use only of the available historical record, possibly corrupted by model error, and the initial state itself. The approach utilizes out-of-sample extension methods  \citep{coifman2006geometric,rabin2012heterogeneous} to define our target observable beyond a designated training period. There are some key differences between this study and that of \citet{zhao2014analog}. First, we consider multivariate data, including SIC, which is intrinsically noisier than SST, given that SIC is a thresholded variable and experiences high variability in the marginal ice zone. We also make use of observational data, which is a relatively short time series compared to GCM model data. In addition to considering kernel eigenfunctions as our target observable, which we will see is well-suited for prediction, we also consider integrated sea ice anomalies, a more challenging observable uninfluenced by the data analysis algorithm. We compare performance of these kernel ensemble analog forecasting techniques to some parametric forecast models (specifically, autoregressive and FEM-VARX models), and find the former to have higher predictive skill than the latter, which often can not even outperform the simple persistence forecast.
	
	\par The rest of the paper is outlined as follows. In Section \ref{sec:methods}, we discuss the mathematical methods used to perform kernel ensemble analog forecasting predictions, and also discuss alternative parametric forecasting methods. In Section \ref{sec:datasets} we describe the data sets used for our experiments, and in Section \ref{sec:results} we present our results. Discussion and concluding remarks are in Section \ref{sec:discussion}.

\section{Methods}
\label{sec:methods}
	\par Our forecasting approach is motivated by using kernels, a pairwise measure of similarity for the state vectors of interest. A simple example of such an object is a Gaussian kernel:
	\begin{equation}
	w(x_i,x_j) = e^{-\frac{\|x_i - x_j\|^2}{\sigma_0}},
	\label{eq:gaussian_kernel}
	\end{equation}
	where $x$ is a state variable of interest, and $\sigma_0$ is a scale parameter that controls the locality of the kernel. The kernels are then used to generate a weighted ensemble of analog forecasts, making use of an available historical record, rather than relying on a single analog. To extend the observable we wish to predict beyond the historical period into the future, we make use of out-of-sample extension techniques. These techniques allow one to extend a function (observable) to a new point $y$ by looking at values of the function at known points $x_i$ close to $y$, where closeness will be determined by the kernels. An important part of any prediction problem is to choose a target observable that is both physically meaningful and exhibits high predictability. As demonstrated in \citet{zhao2014analog}, defining a kernel on the desired space naturally leads to a preferred class of observables that exhibit time-scale separation and good predictability.
	
	\subsection{Nonlinear Laplacian spectral analysis}
	\label{subsec:nlsa}
	\par Nonlinear dynamical systems generally give rise to datasets with low-dimensional nonlinear geometric structures (e.g., attractors). We therefore turn to a data analysis technique, NLSA, that allows us to extract spatio-temporal patterns from data from a high-dimensional non-linear dynamical system, such as a coupled global climate system \citep{giannakis2012comparing, giannakis2012nonlinear, giannakis2013nonlinear, giannakis2014data}. The standard NLSA algorithm is a nonlinear manifold generalization of singular spectrum analysis (SSA) \citep{ghil2002advanced}, where the covariance operator is replaced by a discrete Laplace-Beltrami operator to account for non-linear geometry on the underlying data manifold. The eigenfunctions of this operator then form a convenient orthonormal basis on the data manifold. A key advantage of NLSA is there is no pre-processing of the data needed, such as band-pass filtering or seasonal partitioning.
	
		\subsubsection{Time-lagged embedding}
		\label{subsec:embed}
		\par An important first step in the NLSA algorithm is to perform a time-lagged embedding of the spatio-temporal data as a method of inducing time-scale separation in the extracted modes. An analog forecast method is driven by the initial data, and to incorporate some of the system's dynamics into account in selecting an analog, instead of using a snapshot in time of the state as our initial condition, time-lagged embedding helps make the data more Markovian. 
		\par Let $z(t_i) \in\mathbb{R}^d$ be a time series sampled uniformly with time step $\delta t$, on a grid of size $d$, with $i=1,\ldots, N$ samples. We construct a lag-embedding of the data set with a window of length $q$, and consider the lag-embedded time series 
		\begin{equation*}
		x(t_i) = \left(z(t_i), z(t_{i-1}), \ldots, z(t_{i-(q-1)})\right)\in\mathbb{R}^{dq}.
		\end{equation*}
		The data now lies in $\mathbb{R}^m$, with $m=dq$ the dimension of the lagged embedded space, and $n=N-q+1$ number of samples in lagged embedded space (also called Takens embedding space, or delay coordinate space). It has been shown that time-lagged embedding recovers the topology of the attractor of the underlying dynamical system that has been lost through partial observations \citep{takens1981detecting, sauer1991embedology}. In particular, the embedding affects the non-linear geometry of the underlying data manifold $\mathcal{M}$, in such a way to allow for dynamically stable patterns with time scale separation \citep{berry2013time}, a desirable property that will lead to observables with high predictability.
		
		\subsubsection{Discrete Laplacian}
		\label{subsubsec:discrete_laplacian}
		\par The next step is to define a kernel on the data manifold $\mathcal{M}$. Rather than using a simple Gaussian as in Equation \eqref{eq:gaussian_kernel}, the NLSA kernel makes use of phase velocities $\xi_i = \|x(t_i) - x(t_{i-1})\|$, which forms a vector field on the data manifold and provides additional important dynamic information. The NLSA kernel we use is
		\begin{equation}
		K\left(x(t_i),x(t_j)\right) = \exp\left(-\frac{\|x(t_i) - x(t_j)\|^2}{\epsilon\xi_i\xi_j}\right).
		\label{eq:nlsa_kernel}
		\end{equation}
		With this kernel $K$ and associated matrix $K_{ij} = K\left(x(t_i),x(t_j)\right)$, we solve the Laplacian eigenvalue problem to acquire an eigenfunction basis $\{\phi_i\}$ on the data-manifold $\mathcal{M}$. To do this, we construct the discrete graph Laplacian by following the diffusion maps approach of \citet{coifman2006diffusion}, and forming the following matrices:
		\begin{equation*}
		Q_i = \sum_{j=1}^nK_{ij}, \qquad \tilde{K}_{ij} = \frac{K_{ij}}{Q_i^\alpha Q_j^\alpha},
		\end{equation*}
		\begin{equation*}
		D_i = \sum_{j=1}^n\tilde{K}_{ij}, \qquad P_{ij} = \frac{\tilde{K}_{ij}}{D_i}, \qquad L_{ij} = I - P_{ij}.
		\end{equation*}		
		Here $\alpha$ is a real parameter, typically with value 0, $\tfrac{1}{2}$, 1. We note in the large data limit, as $n\to\infty$ and $\epsilon\to0$, this discrete Laplacian converges to the Laplace-Beltrami operator on $\mathcal{M}$ for a Riemannian metric that depends on the kernel \citep{coifman2006diffusion}. We can therefore think of the kernel as biasing the geometry of the data to reveal a class of features, and the NLSA kernel does this in such a way as to extract dynamically stable modes with time scale separation. We then solve the eigenvalue problem
		\begin{equation}
		L\phi_i = \lambda_i\phi_i.
		\label{eq:eig}
		\end{equation}
		The resulting Laplacian eigenfunctions $\phi_i = (\phi_{1i},\ldots,\phi_{ni})^T$ form an orthonormal basis on the data manifold $\mathcal{M}$ with respect to the weighted inner product:
		\begin{equation*}
		\langle\phi_i,\phi_j\rangle = \sum_{k=1}^nD_k\phi_{ki}\phi_{kj} = \delta_{ij}.
		\label{eq:eig}
		\end{equation*}
		\par As well as forming a convenient orthonormal basis, the eigenfunctions $\phi_i$ give us a natural class of observables with good time-scale separation and high predictability. These eigenfunctions $\phi_i$ are time series, nonlinear analogs to principal components, and can be used to recreate spatio-temporal modes, similar to extended empirical orthogonal functions (EOFs) \citep{giannakis2012limits, giannakis2012nonlinear, giannakis2013nonlinear}. However, unlike EOFs, the eigenfunctions $\phi_i$ do not measure variance, but rather measure oscillations or roughness in the abstract space $\mathcal{M}$, the underlying data manifold. The eigenvalues $\lambda_i$ measure the Dirichlet energy of the corresponding eigenfunctions, which has the interpretation of being squared wave numbers on this manifold \citep{giannakis2014data}. We now use the leading low-frequency kernel eigenfunctions as our target observables for prediction.
	
		\subsubsection{Multiple components}
		\par As described in \citet{bushuk2014reemergence}, the above NLSA algorithm can be modified to incorporate more than one time series. Let $z^{(1)}$, $z^{(2)}$ be two signals, sampled uniformly with time step $\delta t$, on (possibly different) $d_1, d_2$ grid points. After lag-embedding each variable with embedding windows $q_1, q_2$ to its appropriate embedding space, so $z^{(1)}(t_i)\in\mathbb{R}^{d_1}\mapsto x^{(1)}(t_i)\in\mathbb{R}^{d_1q_1}$ and $z^{(2)}(t_i)\in\mathbb{R}^{d_2}\mapsto x^{(2)}(t_i)\in\mathbb{R}^{d_2q_2}$, we construct the kernel function $K$ by scaling the physical variables $x^{(1)}, x^{(2)}$ to be dimensionless by
		\begin{equation}
		K_{ij} = \exp\left(-\frac{\|x^{(1)}(t_i) - x^{(1)}(t_j)\|^2}{\epsilon\xi^{(1)}_i\xi^{(1)}_j} - \frac{\|x^{(2)}(t_i) - x^{(2)}(t_j)\|^2}{\epsilon\xi^{(2)}_i\xi^{(2)}_j}\right).
		\label{eq:kernel_multivariate}
		\end{equation}
		\par This can then be extended to any number of variables, regardless of physical units, and allows for analysis in coupled systems, such as the ocean and sea ice components of a climate model. An alternative approach to Equation \eqref{eq:kernel_multivariate} in extending the NLSA kernel to multiple variables with different physical units is to first normalize each variable to unit variance. However a drawback of this approach is that information about relative variability is lost as the ratios of the variances are set equal, rather than letting the dynamics of the system control the variance ratios of the different variables \citep{bushuk2014reemergence} as in Equation \eqref{eq:kernel_multivariate}.

	\subsection{Out-of-sample extension}
	\label{subsec:OSE}
	\par Now that we have established our class of target observables, namely the eigenfunctions $\phi_i$ in Equation \eqref{eq:eig}, we need a method for extending these observables into the future to form our predictions, for which we draw upon out-of-sample extension techniques. To be precise, let $f$ be a function defined on a set $M = \{x_1,\ldots,x_n\}, x_i\in\mathbb{R}^m$; $f$ may be vector-valued, but in our case is scalar. We wish to make a prediction of $f$ by extending the function to be defined on a point outside of the training set $M$, by performing an out-of-sample extension, which we call $\bar{f}$. There are some desireable qualities we wish to have in such an extension, namely that $\bar{f}$ is in some way well-behaved and smooth on our space, and is in consistent as the number of in-samples increases. Below we discuss two such methods of out-of-sample extension, the geometric harmonics, based on the Nystr{\"o}m method, and Laplacian pyramids.
	 
		\subsubsection{Geometric harmonics}
		\label{subsubsec:geometric}
		\par The first approach for out-of-sample extension is based on the Nystr{\"o}m method \citep{nystrom1930praktische}, recently adapted to machine learning applications \citep{coifman2006geometric}, and is based on representing a function $f$ on $M$ in terms of an eigenfunction basis obtained from a spectral decomposition of a kernel. While we use the NLSA kernel (Equation \ref{eq:nlsa_kernel}), other kernels could be used, another natural choice being a Gaussian kernel. For a general kernel $w:\mathbb{R}^m\times\mathbb{R}^m\to\mathbb{R}$, consider its row-sum normalized counterpart:
	\begin{equation*}
		W(y_i,x_j) = \frac{w(y_i,x_j)}{\sum_jw(y_i,x_j)}.
		\label{eq:rowsum}
	\end{equation*}
	 These kernels have the convenient interpretation of forming discrete probability distributions in the second argument, dependent on the first argument, so $W(y,x) = p_y(x)$, which will be a useful perspective later in our analog forecasting in Section \ref{subsec:KEAF}. We then solve the eigenvalue problem
		\begin{equation*}
		\lambda_l\varphi_l(x_i) = \sum_{j=1}^nW(x_i,x_j)\varphi_l(x_j).
		\end{equation*}
		We note the spectral decomposition of $W$ yields a set of real eigenvalues $\lambda_l$ and an orthonormal set of eigenfunctions $\varphi_l$ that form a basis for $L^2(M)$ \citep{coifman2006geometric}, and we can thus represent our function $f$ in terms of this basis:
		\begin{equation}
		f(x_i) = \sum_{j=1}^n\langle\varphi_j,f\rangle\varphi_j(x_i).
		\label{eq:f_decomp}
		\end{equation}
		Let $y\in\mathbb{R}^m, y\notin M$ be an out-of-sample, or test data, point, to which we wish to extend the function $f$. If $\lambda_l\neq0$, the eigenfunction $\varphi_l$ can be extended to any $y\in\mathbb{R}^m$ by
		\begin{equation}
		\bar{\varphi}_l(y) = \frac{1}{\lambda_l}\sum_{j=1}^nW(y,x_j)\varphi_l(x_j).
		\label{eq:ose_eig}
		\end{equation}
		This definition ensures that the out-of-sample extension is consistent when restricted to $M$, meaning $\bar{\varphi}_l(y) = \varphi_l(y)$ for $y\in M$. Combining Equations \eqref{eq:f_decomp} and \eqref{eq:ose_eig} allows $f$ to be assigned for any $y\in\mathbb{R}^m$ by evaluating the eigenfunctions $\varphi_l$ at $y$ and using the projection of $f$ onto these eigenfunctions as weights:
		\begin{equation}
		\bar{f}(y) = \sum_{j=1}^n\langle\varphi_j,f\rangle\bar{\varphi}_j(y).
		\label{eq:ose_nystrom}
		\end{equation}
		Equation \eqref{eq:ose_nystrom} is called the Nystr{\"o}m extension, and in \citet{coifman2006geometric} the extended eigenfunctions in Equation \eqref{eq:ose_eig} are called geometric harmonics. We note this scheme becomes ill conditioned since $\lambda_l\to0$ as $l\to\infty$ \citep{coifman2006geometric}, so in practice there is a truncation of the sum at some level $l$, usually determined by the decay of the eigenvalues. With the interpretation of the eigenvalues $\lambda$ as wavenumbers on the underlying data manifold, this truncation represents removing features from the data that are highly oscillatory in this space.
				
		\subsubsection{Laplacian pyramid}
		\label{subsubsec:LP}
		\par The geometric harmonics method is well suited for observables that have a tight bandwidth in the eigenfunction basis (particularly the eigenfunctions themselves), but for observables that may require high levels of eigenfunctions in their representation, the above mentioned ill-conditioning may hamper this method. An alternative to geometric harmonics is the Laplacian pyramid \citep{rabin2012heterogeneous}, which invokes a multiscale decomposition of the original function $f$ in its out-of-sample extension approach.
		\par A family of kernels defined at different scales is needed, and for clarity of exposition we will use a family of Gaussian kernels $w_l$ (and their row-normalized counterparts $W_l$) at scales $l$:
		\begin{equation}
		w_l(x_i,x_j) = e^{-\frac{\|x_i - x_j\|^2}{\sigma_0/2^l}},
		\label{eq:LPkernels}
		\end{equation}
		That is, $l=0$ represents the widest kernel width, and increasing $l$ gives finer kernels resolving more localized structures.
		\par For a function $f:M\to\mathbb{R}$, the Laplacian pyramid representation of $f$ approximates $f$ in a multiscale manner by $f\approx s_0 + s_1 + s_2 + \cdots$, where the first level $s_0$ is defined by
		\begin{equation*}
		s_0(x_k) = \sum_{i=1}^nW_0(x_i,x_k)f(x_i),
		\end{equation*}
		and we then evaluate the difference:
		\begin{equation*}
		d_1 = f - s_0.
		\end{equation*}
		We then iteratively define the $l$th level decomposition $s_l$:
		\begin{equation*}
		s_l(x_k) = \sum_{i=1}^n W_l(x_i,x_k)d_l(x_i), \qquad d_l = f - \sum_{i=0}^{l-1}s_i.
		\end{equation*}
		Iteration is continued until some prescribed error tolerance $\|f - \sum_ks_k\|< \varepsilon$ is met.
		\par Next we extend $f$ to a new point $y\in\mathbb{R}^m, y\notin D$ by
		\begin{equation*}
		\bar{s}_0(y) = \sum_{i=1}^n W_0(x_i,y)f(x_i),\qquad \bar{s}_l(y) = \sum_{i=1}^n W_l(x_i,y)d_l(x_i),
		\end{equation*}
		for $l\geq1$, and assign $f$ the value
		\begin{equation}
		\bar{f}(y) = \sum_k \bar{s}_k(y).
		\label{eq:LP}
		\end{equation}
		That is, we have formed a multiscale representation of $f$ using weighted averages of $f$ for nearby inputs, where the weights are given by the scale of the kernel function. Since the kernel function can accept any inputs from $\mathbb{R}^m$, we can define these weights for other points outside $M$, and thus define $f$ by using weighted values of $f$ on $M$ (known), where now the weights are given by the proximity of the out-of-sample $y\notin M$ to input points $x_i\in M$. The parameter choices of the initial scale $\sigma_0$ and error tolerance $\varepsilon$ set the scale and cut-off of the dyadic decomposition of $f$ in the Laplacian pyramid scheme. We choose $\sigma_0$ to be the median of the pairwise distances of our training data, and the error tolerance $\varepsilon$ to be scaled by the norm of the observable over the training data, for example $10^{-6}\|f\|$. In our applications below, we use a multiscale family of NLSA kernels based on Equation \eqref{eq:kernel_multivariate} rather than the above family of Gaussian kernels.
		
	\subsection{Kernel ensemble analog forecasting}
	\label{subsec:KEAF}
	\par The core idea of traditional analog forecasting is to identify a suitable analog to one's current initial state from a historical record, and then make a prediction based on the trajectory of that analog in the historical record. The analog forecasting approach laid out in \citet{zhao2014analog}, which we use here, varies in a few important regards. First, the initial system state, as well as the historical record (training data), is in Takens embedding space, so that an analog is not determined by a snapshot in time alone, but the current state with some history (a `video'). Second, rather than using Euclidean distance, as in traditional analog forecasting, the distances we use are based on a defined kernel function, which reflects a non-Euclidean geometry on the underlying data manifold. The choice of geometry is influenced by the lagged embedding and choice of kernel, which we have done in a way that gives us time scale separation in the resulting eigenfunctions $\phi$, yielding high predictability. Third, rather than identify and use a single analog in the historical record, weighted sums of analogs are used, with weights determined by the kernel function.
	\par For this last point, it is useful to view analog forecasting in the context of an expectation over an empirical probability distribution. For example, a traditional analog forecast of a new initial state $y$ at some time $\tau$ in the future, based on selected analog $x_j = x(t_j)$, can be written as
	\begin{equation*}
	f(y,\tau)= \mathbb{E}_{p_y}S_\tau f = \sum_{i=1}^n p_y(x_i)f\left(x(t_i+\tau)\right) = f\left(x(t_j + \tau)\right),  
	\end{equation*}
	where $p_y = \delta_{ij}$ is the Dirac delta function, and $S_\tau f(x_i) = f(x(t_i+\tau))$ is the operator that shifts the timestamp of $x_i$ by $\tau$. To make use of more than one analog and move to an ensemble, we let $p_y$ be a more general discrete empirical distribution, dependent on the initial condition $y$, with probabilities (weights) determined by our kernel function. Writing in this way, we simply need to define the empirical probability distribution $p_y$ to form our ensemble analog forecast.
	\par For the geometric harmonics method, we projected $f$ onto an eigenfunction basis $\phi_i$ (truncated at some level $l$) and performed an out-of-sample extension on each eigenfunction basis function, i.e.
	\begin{equation*}
	\bar{f}(y) = \sum_{i=1}^l\langle\phi_i,f\rangle\bar{\phi}_i(y).
	\end{equation*}
	Or, to write this in terms of an expectation over an empirical probability measure:
	\begin{equation*}
	\bar{f}(y) = \mathbb{E}_{p_y}f = \sum_{i=1}^l\langle\phi_i,f\rangle \mathbb{E}_{p_y}\phi_i(y),
	\end{equation*}
	where
	\begin{equation*}
	\mathbb{E}_{p_y}\phi_i = \frac{1}{\lambda_i}\sum_{j=1}^nW(y,x_j)\phi_i(x(t_j)).
	\end{equation*}
	We then define our prediction for lead time $\tau$ via geometric harmonics by
	\begin{equation*}
	f(y,\tau)  =  \mathbb{E}_{p_y}S_\tau f = \sum_{i=1}^l\langle \phi_i, f\rangle \mathbb{E}_{p_y}S_\tau \phi_i,
	\end{equation*}
	where
	\begin{equation*}
	\mathbb{E}_{p_y}S_\tau \phi_i = \frac{1}{\lambda_i}\sum_{j=1}^nW(y,x_j)\phi_i(x(t_j+\tau)).
	\end{equation*}
	\par Similarly, the $\tau$ shifted ensemble analog forecast via Laplacian pyramids is then
	\begin{equation*}
	f(y,\tau) = \mathbb{E}_{p_{y,0}}S_\tau f + \sum_{i=1}^l\mathbb{E}_{p_{y,i}}S_\tau d_i,
	\end{equation*}
	where $p_{y,i}(x) = W_i(y,x)$ corresponds to the probability distribution from the kernel at scale $i$.
	\par Thus we have a method for forming a weighted ensemble of predictions, that is non-parametric and data-driven through the use of a historical record (training data), which itself has been subject to analysis that reflects the dynamics of the high-dimensional system in the non-linear geometry on the underlying abstract data manifold $\mathcal{M}$, and produces a natural preferred class of observables to target for prediction through the kernel eigenfunctions. 

	\subsection{Autoregressive modeling}
	\label{subsec:ARmodel}
	\par We wish to compare our ensemble analog prediction methods to more traditional parametric methods, namely autoregressive models of the North Pacific variability \citep{frankignoul1977stochastic}, of the form
	\begin{equation*}
	x(t+1) = \mu(t) + A(t)x(t) + \sigma(t)\epsilon(t)
	\label{eq:ARmodel}
	\end{equation*}
	where $x(t)$ is our signal, $\mu(t)$ is the external forcing (possibly 0), $A(t)$ is the autoregressive term, and $\sigma(t)$ is the noise term, each to be estimated from the training data, and $\epsilon(t)$ is a Gaussian process. In the stationary case, the model coefficients $\mu(t) = \mu$, $A(t) = A$, $\sigma(t)=\sigma$ are constant in time, and can be evaluated in an optimal way through ordinary least squares. In a non-stationary case, we will invoke the FEM-VARX framework of \citet{horenko2010identification} by clustering the training data into $K$ clusters, and evaluating coefficients $A_k$, $\sigma_k$, $k=1,\ldots,K$ for each cluster (see \citet{horenko2010clustering, horenko2010identification} for more details on this algorithm). From the algorithm we obtain a cluster identification function $\Gamma(t)$, such that $\Gamma(t_i) = k$ indicates at time $t_i$ the model coefficients are $A(t_i) = A_k$, $\sigma(t_i) = \sigma_k$. In addition to choosing the number of clusters $K$, the method also has a persistence parameter $C$ that governs the number of allowable switches between clusters that must be chosen prior to calculating model coefficients. These parameters are usually chosen to be optimal in the sense of the Akaike Information Criterion (AIC) \citep{horenko2010identification,metzner2012analysis}, an information theoretic based measure for model selection which penalizes overfitting by large number of parameters.
	\par As mentioned earlier, while non-stationary autoregressive models may perform better than stationary models in fitting, an inherent difficulty in a prediction setting is the advancement of the model coefficients $A(t), \sigma(t)$ beyond the training period, which in the above framework amounts to solely advancing the cluster affiliation function $\Gamma(t)$. If we call $\pi_k(t)$ the probability of the model being at cluster $k$ at time $t$, we can view $\Gamma(t)$ as determining a Markov switching process on the cluster member probabilities 
	\begin{equation*}
	\pi(t) = (\pi_1(t),\ldots,\pi_K(t)),
	\end{equation*}
	which over the training period will be 1 in one entry, and 0 elsewhere at any given time. We can estimate the transition probability matrix $T$ of that Markov process by using the optimal cluster affiliation sequence $\Gamma(t)$ from the FEM-VARX framework (here optimal is for the training period). Assuming the Markov hypothesis, we can estimate the stationary probability transition matrix directly from $\Gamma(t)$ by:
	\begin{equation*}
	T_{ij} = \frac{N_{ij}}{\sum_{k=1}^KN_{ik}},
	\end{equation*}
	where $N_{ij}$ is the number of direct transitions from state $i$ to state $j$ \citep{horenko2011nonstationarity, franzke2009systematic}. This estimated transition probability matrix $T$ can be used to model the Markov switching process in the following ways.
		\subsubsection{Generating predictions of cluster affiliation}
		\label{subsubsec:pred_pi}
		\par The first method we employ is to advance the cluster member probabilities $\pi(t)$ using the estimated probability transition matrix $T$ by the deterministic equation \citep{franzke2009systematic}:
		\begin{equation*}
		\pi(t_0 + \tau) = \pi(t_0)T^\tau
		\end{equation*}
		where $\pi(t_0) = (\pi_1(t_0),\pi_2(t_0))$ is the initial cluster affiliation, which is determined by which cluster center the initial point $x(t_0)$ is closest to, and $\pi_i(t_0)$ is either 0 or 1.
		\par The second method we employ is to use the estimated transition matrix $T$ to generate a realization of the Markov switching process $\Gamma_R$, and use this to determine the model cluster at any given time, maintaining strict model affiliation. Thus $\pi_k(t) = 1$ is $\Gamma_R(t) = k$, and 0 otherwise.
		
	\subsection{Error metrics}
	\label{subsec:errormetrics}
	\par To gauge the fidelity of our predictions, we will evaluate the average root-mean-square error (RMSE) and pattern correlation (PC) of our predictions $y$ against the ground truth $x$, where points in our test data set (of length $n'$) are used as intial conditions for our predictions. As a benchmark, we will compare each prediction approach to a simple persistence forecast $y(\tau) = y(0)$. The error metrics are calculated as
	\begin{eqnarray*}
	rms^2(\tau) & = & \frac{1}{n'}\sum_{j=1}^{n'} \left(y(t_j+\tau) - x(t_j+\tau)\right)^2,\\
	pc(\tau) & = & \frac{1}{n'}\sum_{j=1}^{n'}\frac{\left(y(t_j+\tau) - \tilde{y}(\tau)\right)\left(x(t_j+\tau) - \tilde{x}(\tau)\right)}{\sigma_y(\tau)\sigma_x(\tau)},
	\end{eqnarray*}
	where 
	\begin{equation*}
	\tilde{y}(\tau) = \frac{1}{n'}\sum_{j=1}^{n'}y(t_j+\tau), \qquad \tilde{x}(t) = \frac{1}{n'}\sum_{j=1}^{n'} x(t_j+\tau),
	\end{equation*} 
	\begin{eqnarray*}
	\sigma^2_y(\tau) & = & \frac{1}{n'}\sum_{j=1}^{n'}(y(t_j+\tau) - \tilde{y}(\tau))^2, \\
	\sigma^2_x(\tau) & = & \frac{1}{n'}\sum_{j=1}^{n'}(x(t_j+\tau) - \tilde{x}(\tau))^2.
	\end{eqnarray*}
	An important note is that for data-driven observables, such as NLSA eigenfunctions or EOF principal components, there is no underlying ground truth when predicting into the future. As such a ground truth for comparison needs to be defined when evaluating the error metrics, for which one can use the out-of-sample extended function $\bar{f}(y(t_j+\tau))$ as defined in Equations \eqref{eq:ose_nystrom} and \eqref{eq:LP}.

\section{Datasets}
\label{sec:datasets}
	\subsection{CCSM model output}
	\label{subsec:ccsmdata}
	\par We use model data from the Community Climate System Model (CCSM), versions 3 and 4, for monthly SIC and SST data, restricted to the North Pacific, which we define as 20$^\circ$--65$^\circ$N, 120$^\circ$E--110$^\circ$W. CCSM3 model data is used from a 900 year control run (experiment b30.004) \citep{collins2006community}. The sea ice component is the Community Sea Ice Model \citep[CSIM;][]{holland2006influence} and the ocean component is the Parallel Ocean Program (POP), both of which are sampled on the same nominal $1^\circ$ grid. CCSM4 model data is used from a 900 year control run (experiment b40.1850), which uses the Community Ice CodE 4 model for sea ice \citep[CICE4;][]{hunke2008cice} and the POP2 model for the ocean component \citep{smith2010parallel}, also on a common nominal $1^\circ$ grid. Specific differences and improvements between the two model versions can be found in \citet{gent2011community}.
	 \par NLSA was performed on these data sets, both in single and multiple component settings, with the same embedding window of $q_1=q_2=q=24$ months for each variable, kernel scale $\epsilon = 2$, and kernel normalization $\alpha = 0$. A 24 month embedding window was chosen to allow for dynamical memory beyond the seasonal cycle, as is used in other NLSA studies (e.g., \citet{giannakis2012limits, giannakis2013nonlinear, bushuk2014reemergence}). There are 6648 grid points for SST and 3743 for SIC in this region for these models, and with an embedding window of $q=24$, this means our lagged embedded data lies in $\mathbb{R}^m$ with $m=\text{159,552}$ for SST and $m=\text{89,832}$ for SIC. For the purposes of a perfect model experiment, where the same model run is used for both the training and test data, we split the CCSM4 control run into two 400 year sets; years 100 -- 499 for the training data set, and years 500 -- 899 for out-of-sample test points. After embedding, this leaves us with $n = n' = 4777$ samples in each data set. In our model error experiment, we train on 800 years of the CCSM3 control run, and use 800 years of CCSM4 for test data, giving us $n = n' = 9577$ data points.
	 \par In addition to low-frequency kernel eigenfunctions as observables, we also consider North Pacific integrated sea ice extent anomalies as a target for prediction in Section \ref{subsec:siea}. This observable exhibits faster variability, on an intra-annual time-scale, and as such we reduce the embedding window from 24 months to 6 months for our kernel evaluation. Using the CCSM4 model training data, we define monthly anomalies by calculating a climatology $f_c$ of monthly mean sea ice extent. Let $v_j$ be gridpoints in our domain of interest (North Pacific), $c$ the mean SIC, $a$ the grid cell area, and then the sea ice extent anomaly hat will be our target observable is defined as
	\begin{equation}
	f(t_i) = \sum_{j}c(v_j,t_i) a(v_j) - f_c(t_i).
	\label{eq:siea}
	\end{equation}
	 
	 \subsection{Observational data}
	 \label{subsec:obsdata}
	 \par For observational data, we turn to the Met Office Hadley Centre's HadISST data set \citep{rayner2003global}, and use monthly data for SIC and SST, from years 1979--2012, sampled on a 1$^\circ$ latitude-longitude grid. We assign ice covered grid points an SST value of $-1.8^\circ$C, and have removed a trend from the data by calculating a linear trend for each month. There are 4161 spatial grid points for SST, for a lagged embedded dimension of $m=\text{99,864}$, and 3919 grid points for SIC, yielding $m=\text{94,056}$. For direct comparison with this observational data set, the above CCSM model data sets have been interpolated from the native POP grid $1^\circ$ grid to a common 1$^\circ$ latitude-longitude grid. After embedding, we are left with $n'=381$ observation test data points.

\section{Results}
\label{sec:results}

	\subsection{Low-frequency NLSA modes}
	\par
\begin{figure}[ht]
\begin{center}
\includegraphics[width=84mm]{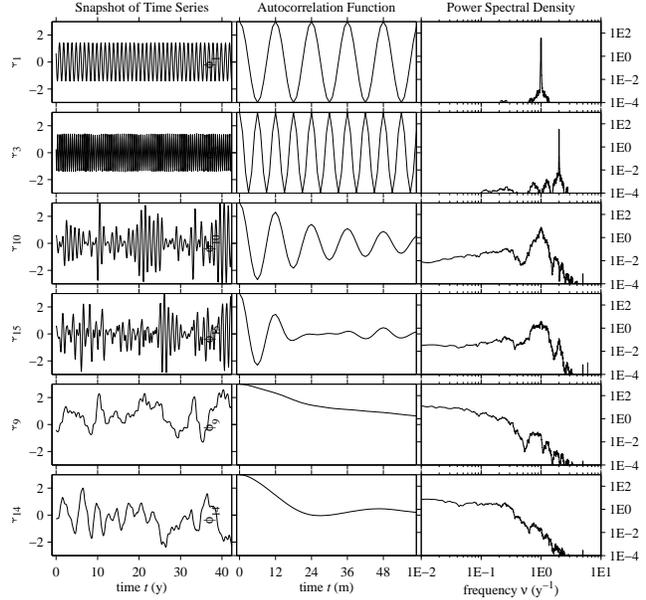}
\caption[NLSA Eigenfunctions]{Select NLSA eigenfunctions for CCSM4 North Pacific SST and SIC data. $\phi_1$, $\phi_3$ are periodic (annual and semi-annual) modes, characterized by a peak in the power spectrum and oscillatory autocorrelation function; $\phi_9$, $\phi_{14}$ are the two leading low-frequency modes, characterized by a red power spectrum and slowly decaying montone autocorrelation function; $\phi_{10}$, $\phi_{15}$ are intermittent modes, characterized by a broad peak in the power spectrum and decaying oscillatory autocorrelation function.}
\label{fig:ccsm4_phi_sstice}
\end{center}
\end{figure}

\begin{figure}[ht]
\begin{center}
\includegraphics[width=84mm]{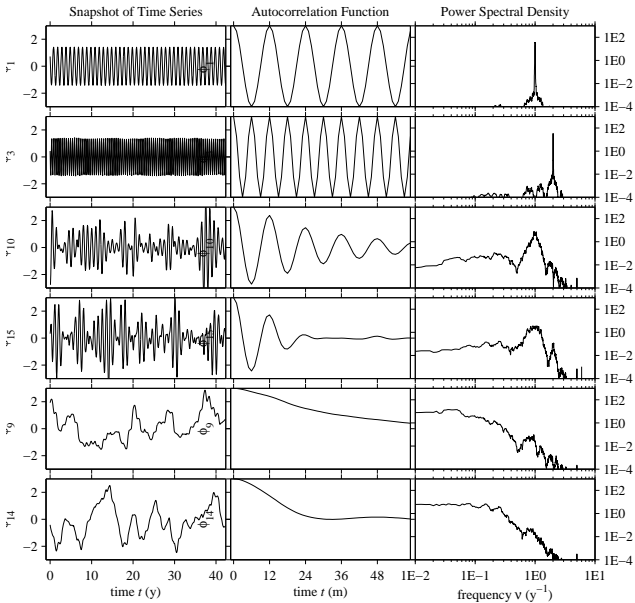}
\caption[NLSA OSE Eigenfunctions]{Out-of-sample extension eigenfunctions, computed via Equation \eqref{eq:ose_nystrom}, for CCSM North Pacific SST and SIC data, associated with the same select (in-sample) eigenfunctions shown in Figure \ref{fig:ccsm4_phi_sstice}.}
\label{fig:ccsm4_osephi_sstice}
\end{center}
\end{figure}

 The eigenfunctions that arise from NLSA typically fall into one of three categories: i) periodic modes, which capture the seasonal cycle and its higher harmonics; ii) low-frequency modes, characterized by a red power spectrum and a slowly decaying autocorrelation function; and iii) intermittent modes, which have the structure of periodic modes modulated with a low-frequency envelope. The intermittent modes are dynamically important \citep{giannakis2012comparing}, shifting between periods of high activity and quiessence, but carry little variance, and are thus typically missed or mixed between modes in classical SSA. Examples of each of these eigenfunctions arising from a CCSM4 data set with SIC and SST variables are shown in Figure \ref{fig:ccsm4_phi_sstice}, for years 100--499 of the pre-industrial control run. The corresponding out-of-sample extension eigenfunctions, defined through the Nystr{\"o}m method, are shown in Figure \ref{fig:ccsm4_osephi_sstice}, and are computed using years 500--899 of the same CCSM4 pre-industrial control run as test (out-of-sample) data. We use the notation $\phi^S_{L_1}$, $\phi^I_{L_1}$, or $\phi^{SI}_{L_1}$, to indicate if the NLSA mode is from SIC, SST, or joint SST and SIC variables, respectively. 
	
	\par We perform our prediction schemes for five year time leads by applying the kernel ensemble analog forecast methods discussed in Section \ref{subsec:KEAF} to the leading two low-frequency modes $\phi^{SI}_{L_1}$, $\phi^{SI}_{L_2}$ from NLSA on North Pacific, shown in Figure \ref{fig:ccsm4_phi_sstice}. The leading low-frequency modes extracted through NLSA can be though of as analogs to the well known PDO and NPGO modes, even in the multivariate setting. We have high correlations between the leading multivariate and univariate low-frequency NLSA modes, with corr $(\phi_{L_1}^{SI},\phi_{L_1}^I) = -0.9907$ for our analog of the NPGO mode, and corr $(\phi_{L_2}^{SI},\phi_{L_1}^S) = 0.8415$ for our analog of the PDO mode.
	\par As a benchmark, we compare against the simple constant persistence forecast $y(t) = y(0)$, which can perform reasonably well given the long decorrelation time of these low-frequency modes, and in fact beats paramteric autoregressive models as we will see below. We define the ground truth itself to be the out-of-sample eigenfunction calculated by Equation \eqref{eq:ose_nystrom}, so by construction, our predictions by the geometric harmonics method are exact at time lag $\tau=0$, whereas predictions using Laplacian pyramids will have reconstruction errors at time lag $\tau=0$.
	
		\subsubsection{Perfect model}
		\par
\begin{figure}[ht]
\begin{center}
\includegraphics[width=39mm]{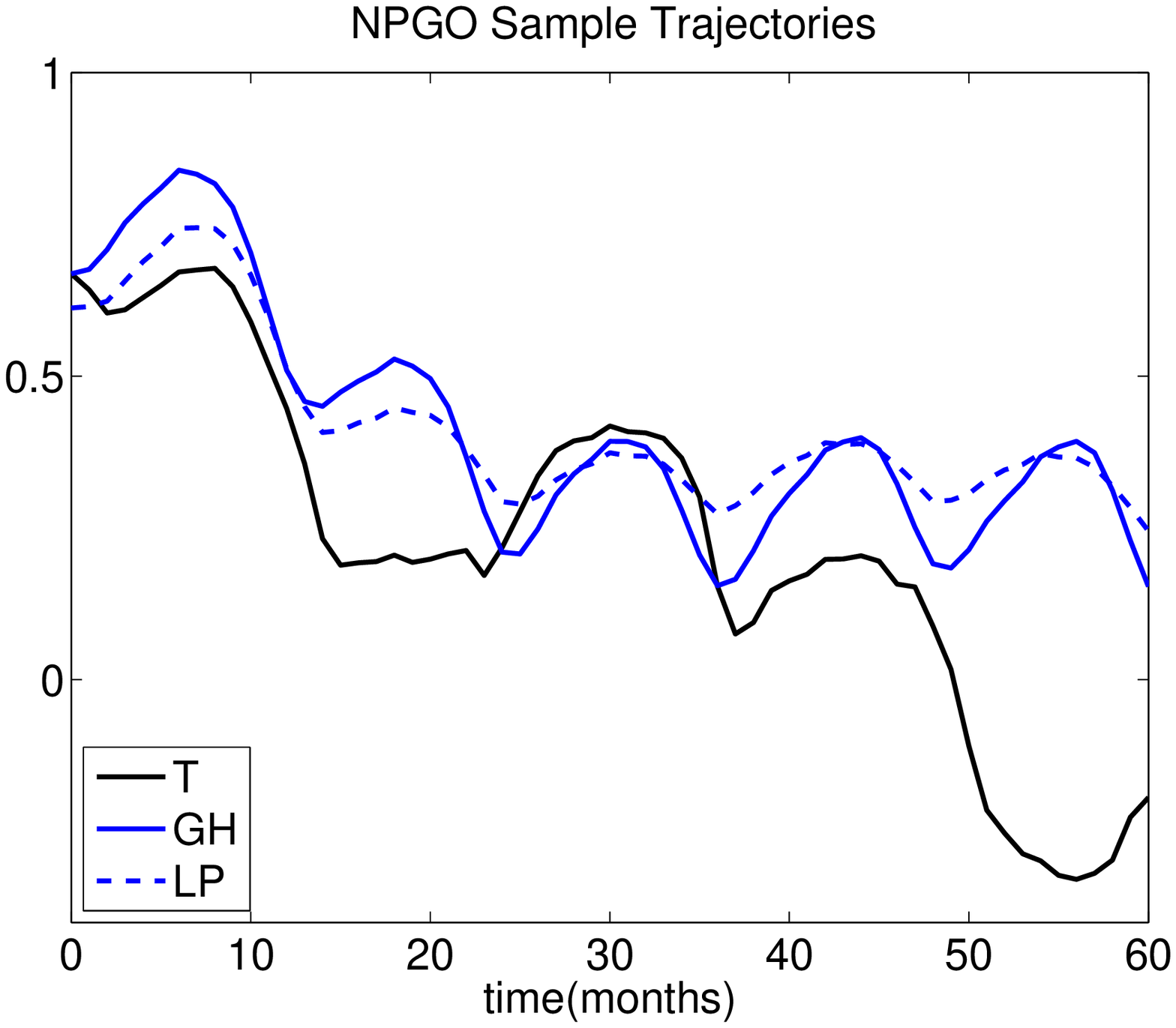}
\includegraphics[width=39mm]{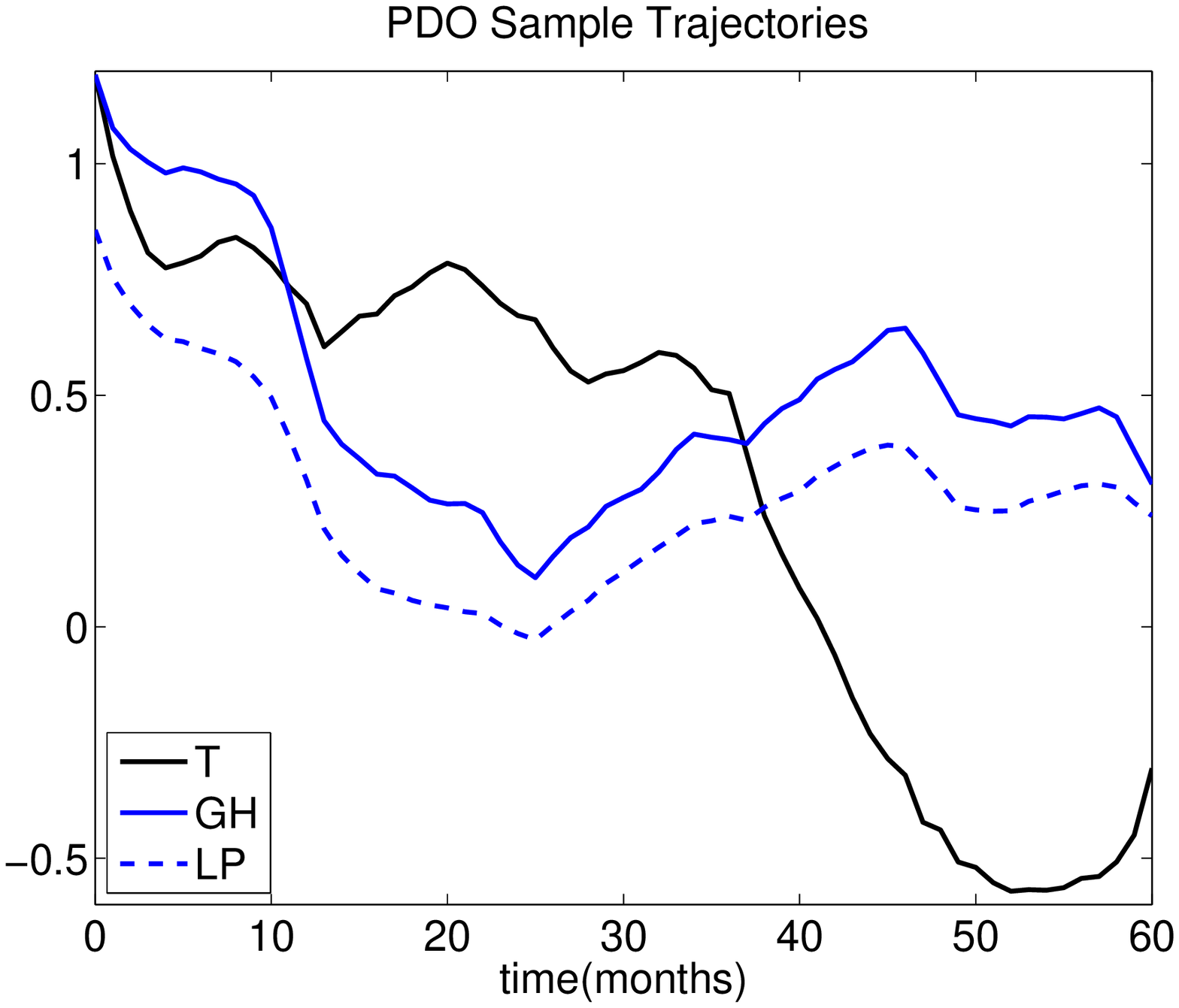}
\caption[KEAF - Perfect Model - Trajectories]{Sample snapshot trajectories of ensemble analog prediction results via geometric harmonics (GH, blue, solid) and Laplacian pyramid (LP, blue, dashed), for the leading two low-frequency modes, in perfect model setting using CCSM4 data. The ground truth (T, black) is itself an out-of-sample extension eigenfunction, as shown in Figure \ref{fig:ccsm4_osephi_sstice}.}
\label{fig:ccsm4_ccsm4_sstice_traj}
\end{center}
\end{figure}

\begin{figure}[ht]
\begin{center}
\includegraphics[width=39mm]{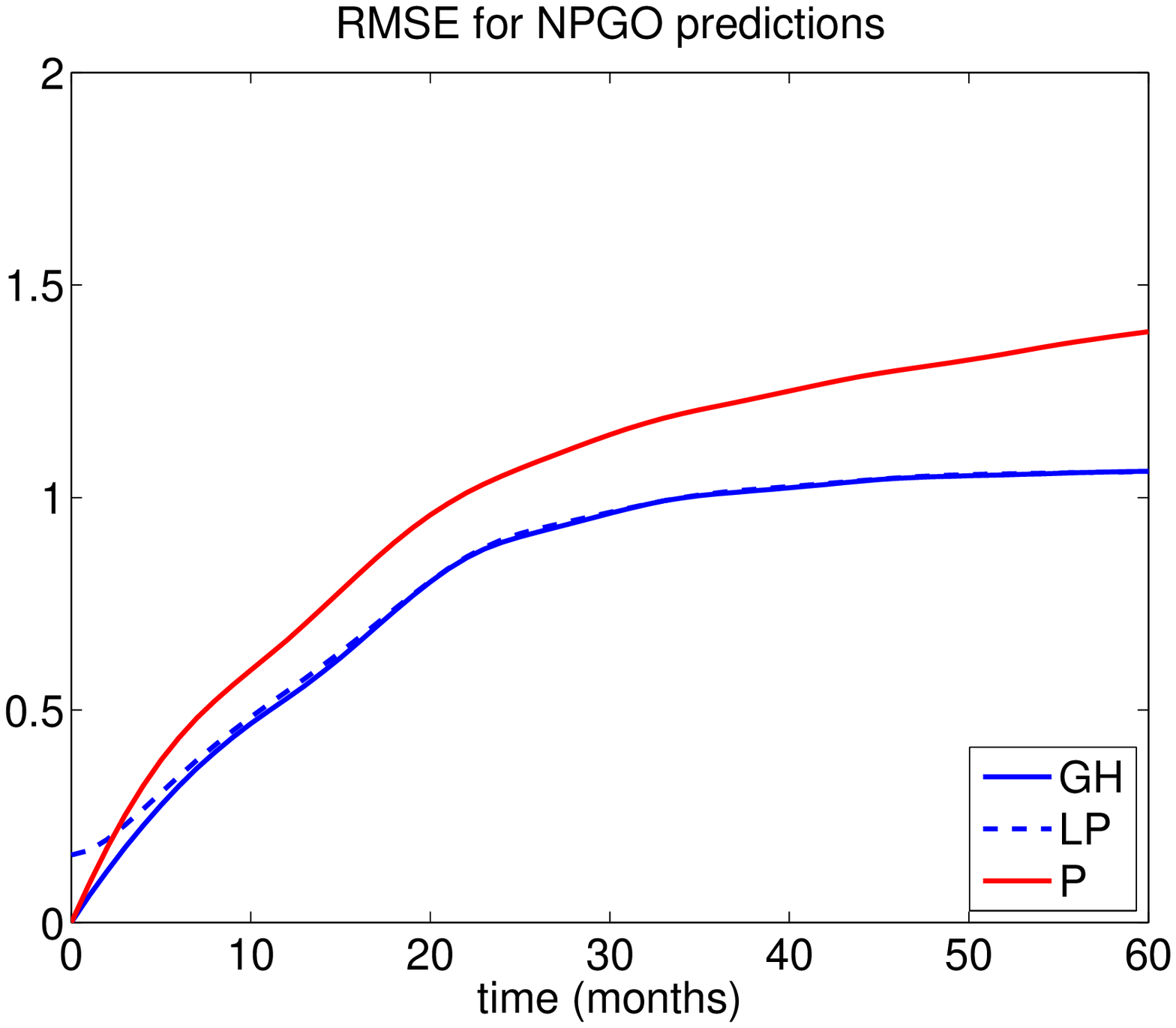}
\includegraphics[width=39mm]{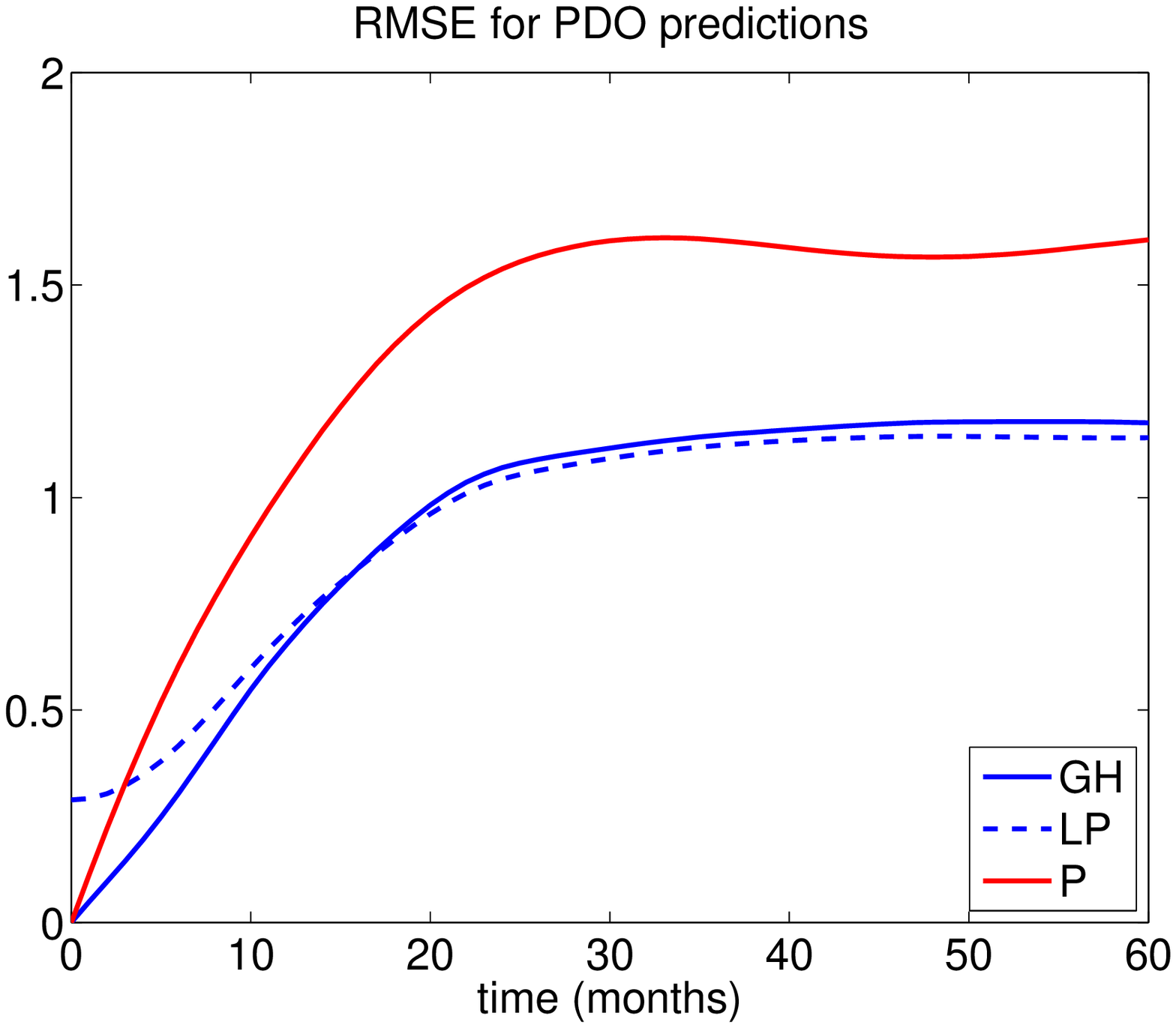}
\includegraphics[width=39mm]{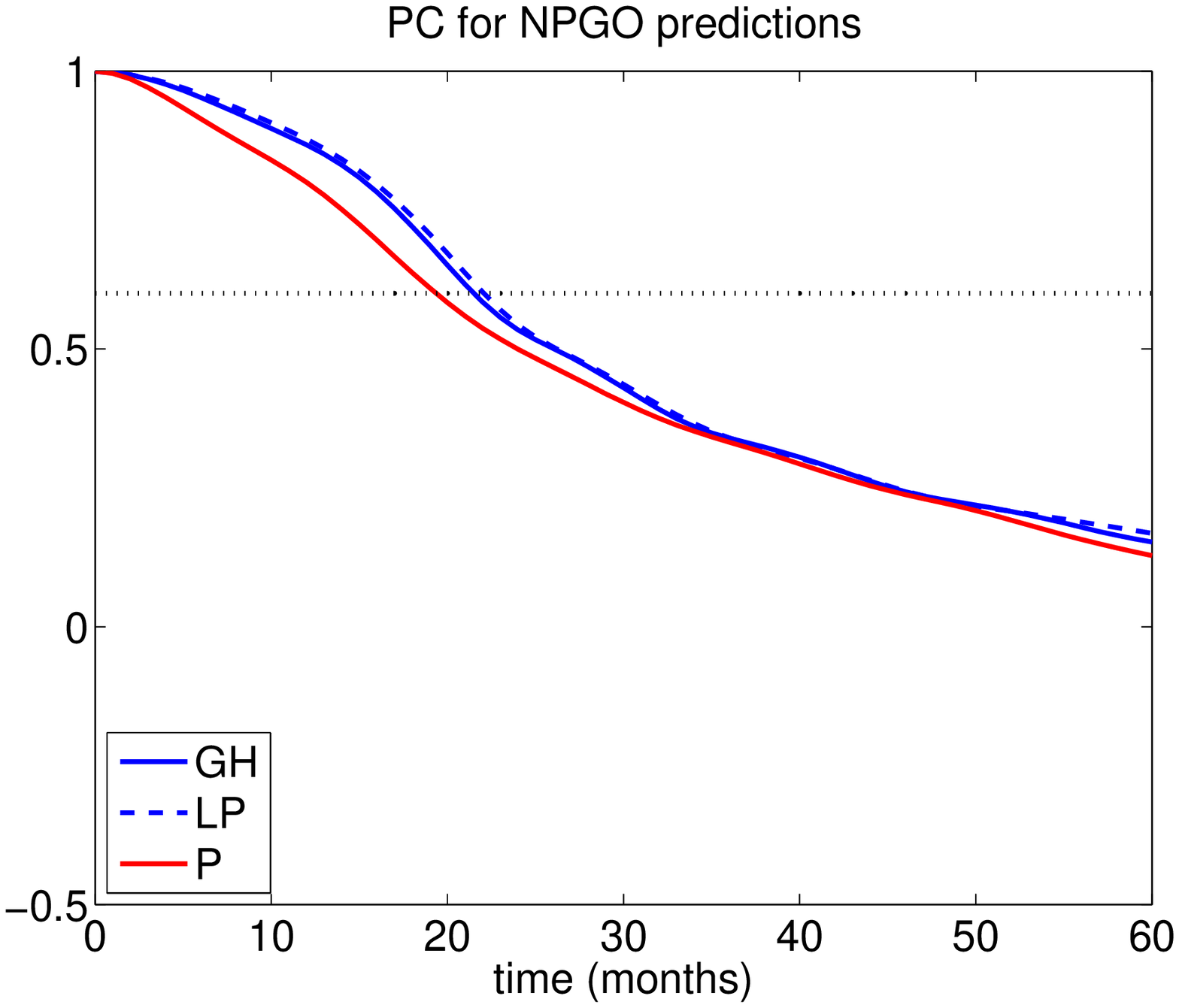}
\includegraphics[width=39mm]{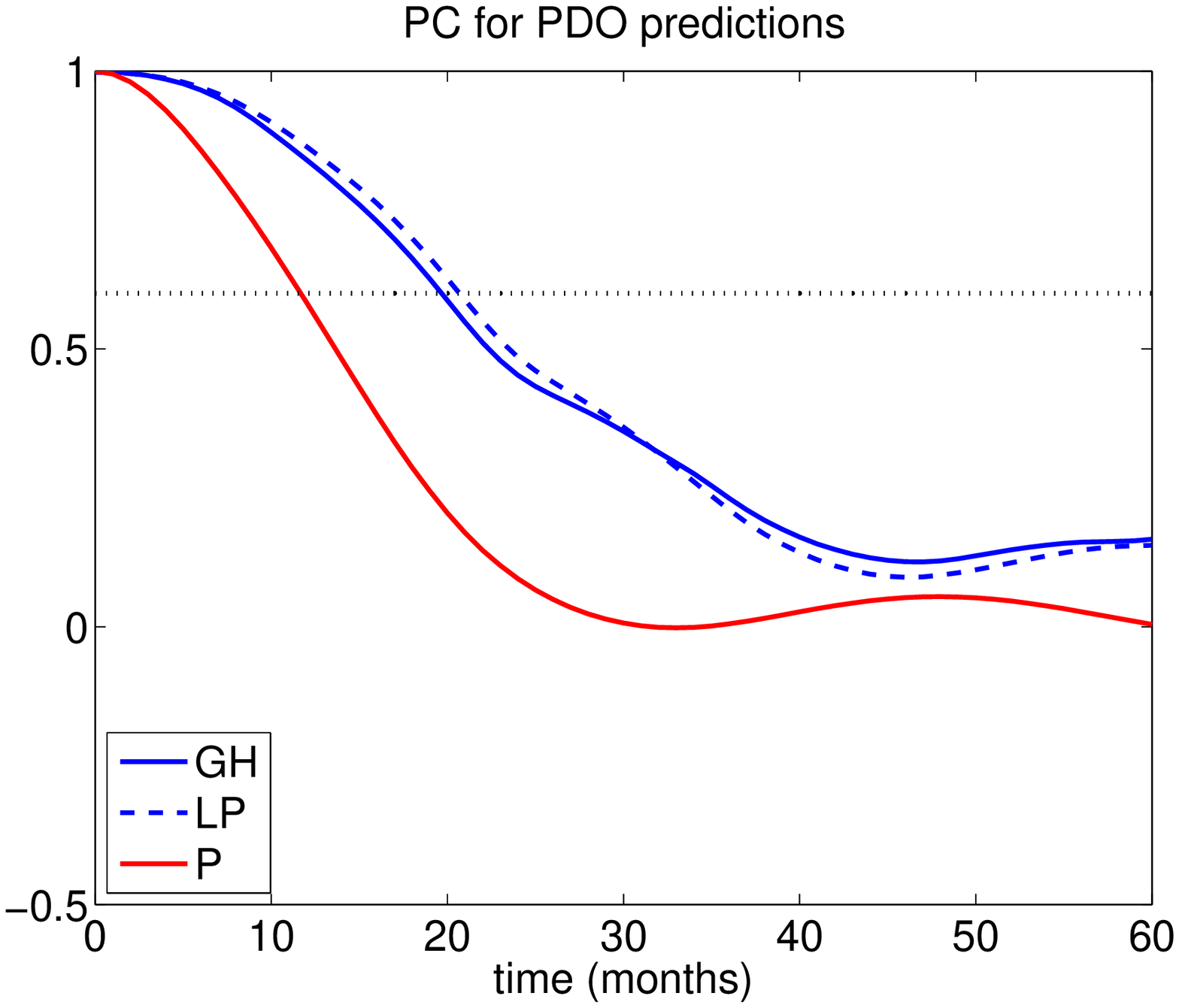}
\caption[KEAF - Perfect Model - Error Metrics]{Kernel ensemble analog prediction results via geometric harmonics and Laplacian pyramid for the leading two low-frequency modes, in perfect model setting using CCSM4 data. Both out-of-sample extension methods outperform the persistence forecast (P) in both error metrics, particularly for in the $\phi^{SI}_{L_2}$ (PDO) mode.}
\label{fig:ccsm4_ccsm4_sstice_error}
\end{center}
\end{figure}

 We first consider the perfect model setting, where the same dynamics generate the training data and test (forecast) data. This should give us a measure of the potential predictability of the methods. Snapshots of sample prediction trajectories along with the associated ground truth out-of-sample eigenfunction are shown in Figure \ref{fig:ccsm4_ccsm4_sstice_traj}. In Figure \ref{fig:ccsm4_ccsm4_sstice_error}, we see that for the leading low-frequency mode $\phi_{L_1}^{SI}$ (our NPGO analog), the ensemble based predictions perform only marginally better than persistence in the PC metric, but have improved RMSE scores over longer timescales. However with the second mode $\phi_{L_2}^{SI}$ (our PDO analog), we see a more noticeable gain in predictive skill with the ensemble analog methods over persistence. If we take 0.6 as a PC threshold \citep{collins2002climate}, below which we no longer consider the model to have predictive skill, we see an increase of about 8 months in predictive skill with the ensemble analog methods over persistence, with skillful forecasts up to 20 months lead time. We note these low-frequency modes extracted from multivariate data exhibit similiar predictability (as measured by when the PC falls below the 0.6 threshold) than their univariate counterparts ($\phi_{L_1}^S$ or $\phi_{L_1}^I$, results not shown).

		\subsubsection{Model error}
		\par
\begin{figure}[ht]
\begin{center}
\includegraphics[width=39mm]{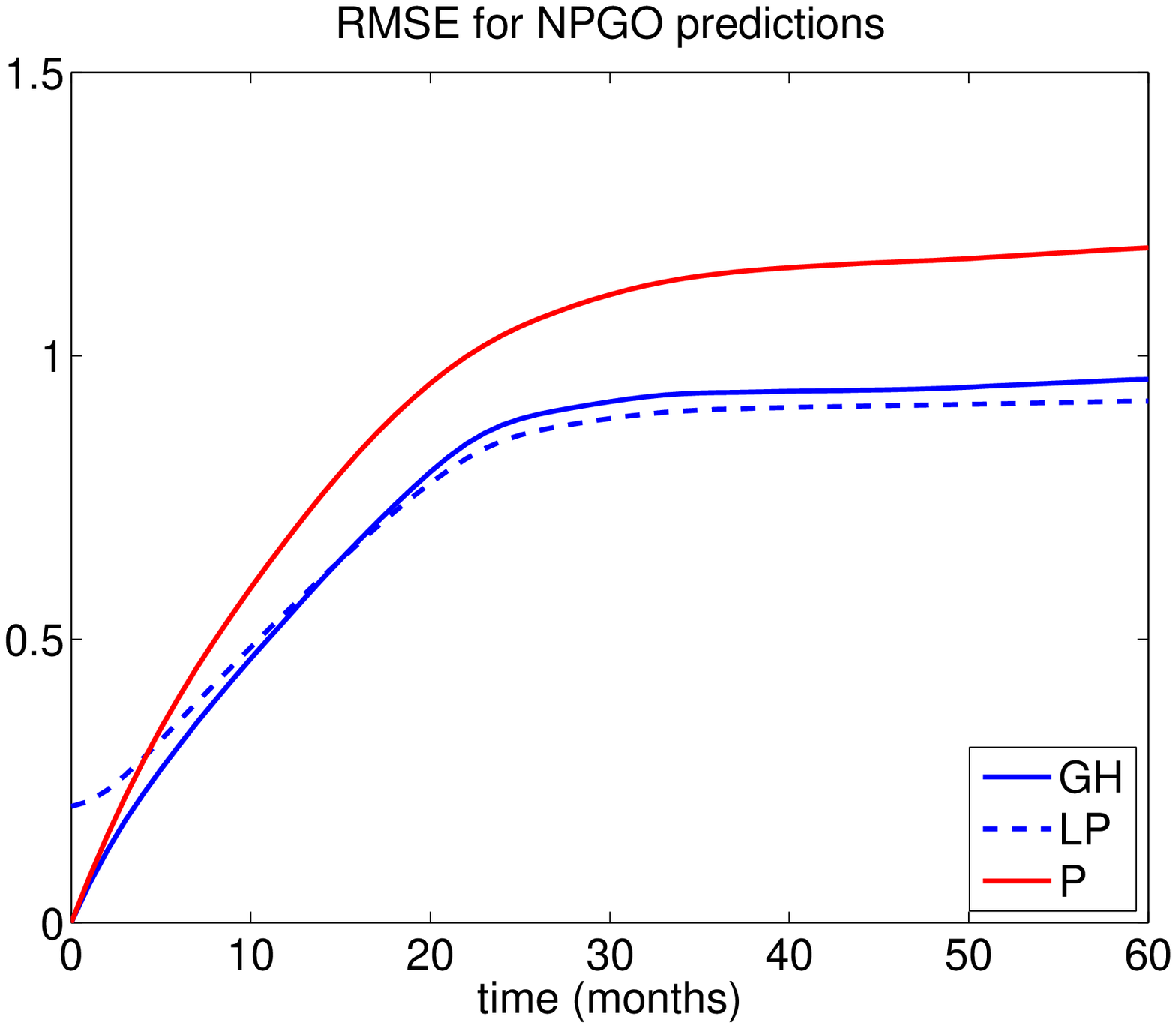}
\includegraphics[width=39mm]{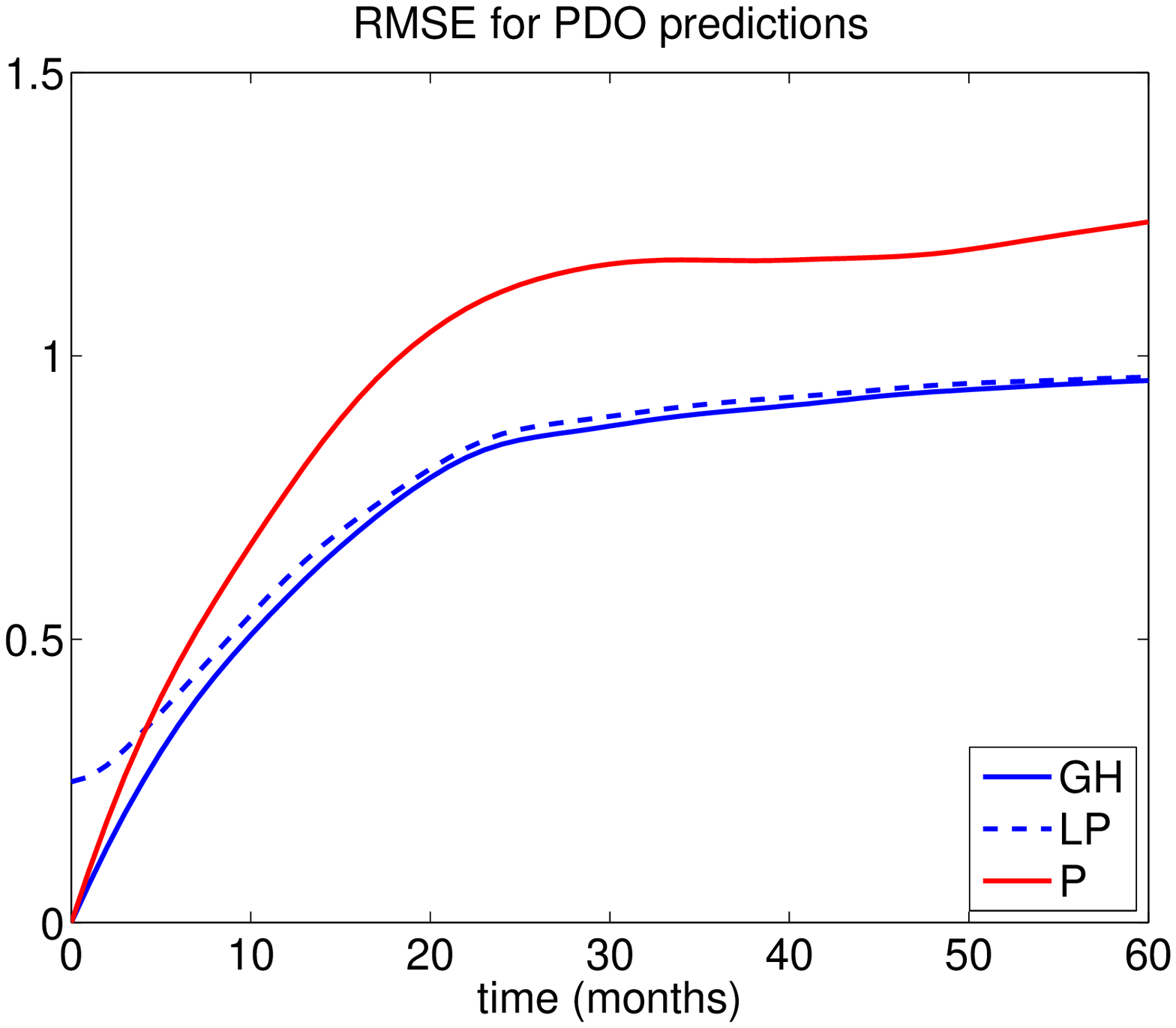}
\includegraphics[width=39mm]{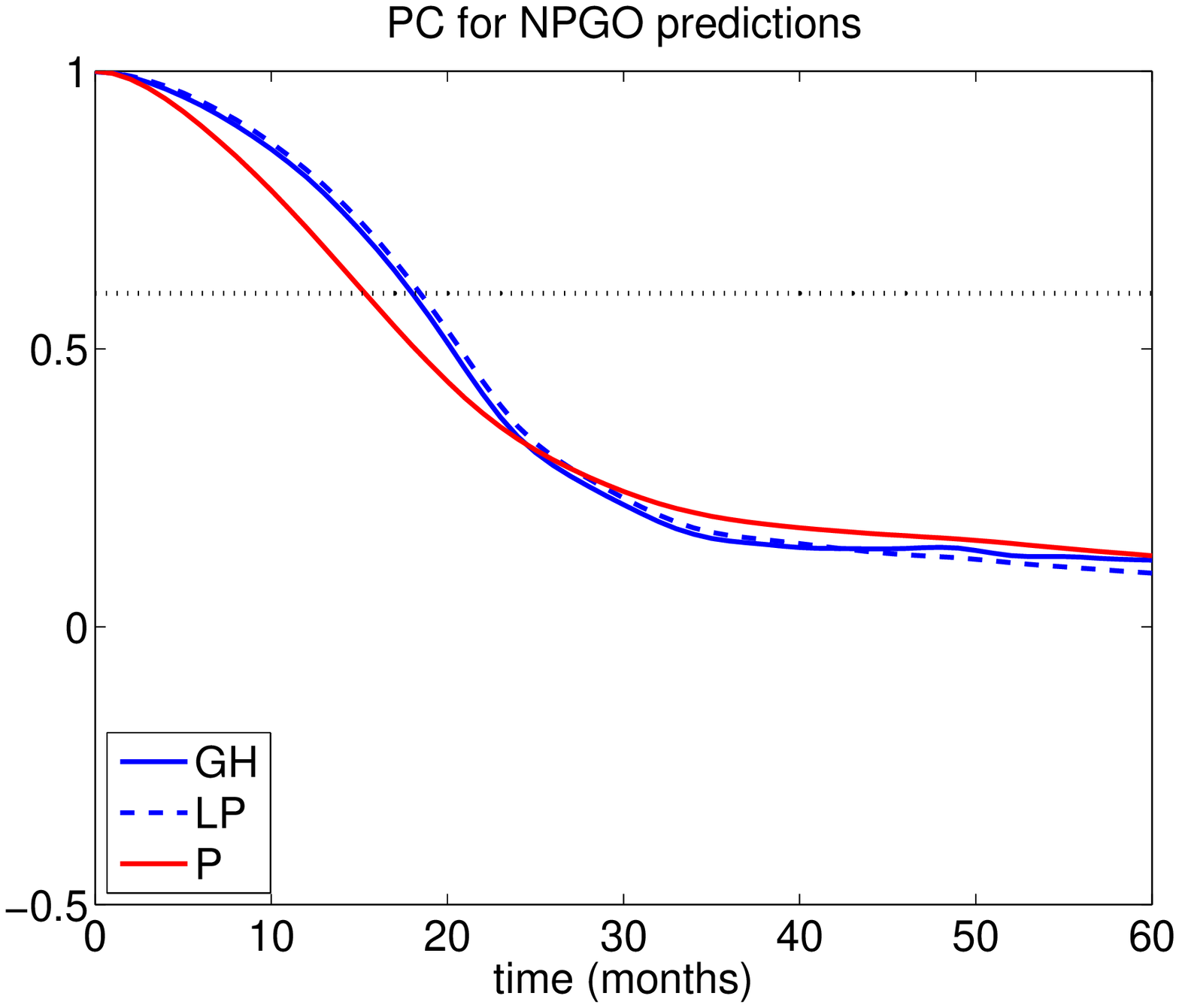}
\includegraphics[width=39mm]{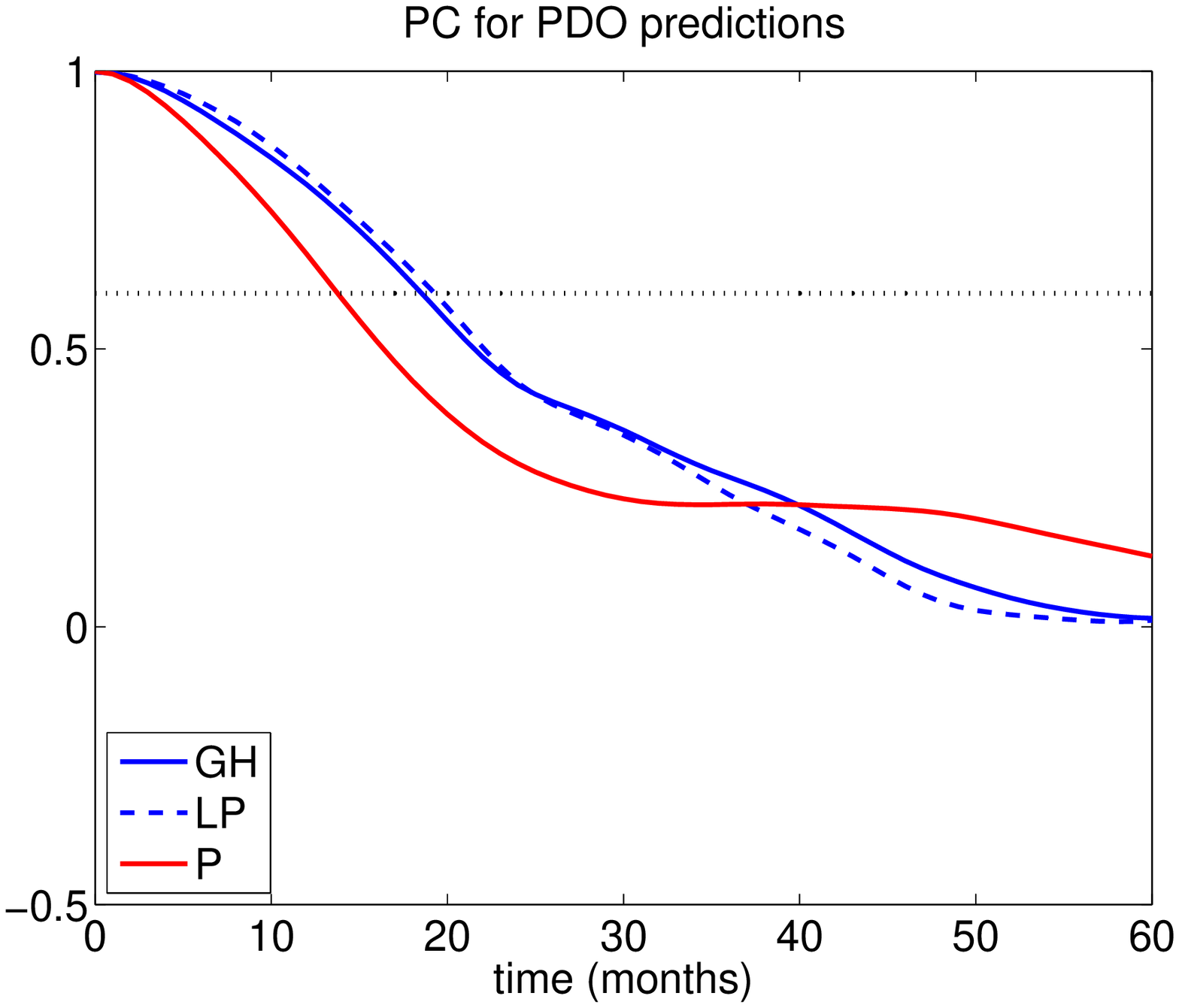}
\caption[KEAF - Model Error - Error Metrics]{Kernel ensemble analog forecasting prediction results for the leading two low-frequency modes, with model error. CCSM3 model data is used for the in-sample data, and CCSM4 model data is used as the out-of-sample data. While the predictions still outperform persistence in the error metrics, there is less gain in predictive skill over as compared to the perfect model case.}
\label{fig:ccsm3_ccsm4_sstice_error}
\end{center}
\end{figure}

 To incorporate model error into our prediction experiment, we train on the CCSM3 model data, serving as our `model', and then use CCSM4 model data as our test data, serving the role of our `nature', and the difference in the fidelity of the two model versions represents general model error. In this experiment, our ground truth is an out-of-sample eigenfunction trained on CCSM3 data, and extended using CCSM4 test data. For the leading $\phi_{L_1}^{SI}$ (NPGO) mode, we see again marginal increased predictive performance in the ensemble analog predictions over persistence at short time scales in PC in Figure \ref{fig:ccsm3_ccsm4_sstice_error}, but at medium to long time scales this improvement has been lost (though after the score has fallen below the 0.6 threshold). This could be due to the increased fidelity of the CCSM4 sea ice model component over the CCSM3 counterpart, where using the less sophisticated model data for training leaves us a bit handicapped in trying to predict the more sophisticated model data \citep{gent2011community,holland2012improved}. The improvement in the predictive skill of $\phi_{L_2}^{SI}$ (PDO) mode over persistence is less pronounced in the presence of model error than it was in the perfect model case shown in Figure \ref{fig:ccsm4_ccsm4_sstice_error}. Nevertheless, the kernel ensemble analog forecasts still provide a substantial improvement of skill compared to the persistence forecast, extending the PC = 0.6 threshold to 20 months.

		\par
\begin{figure}[ht]
\begin{center}
\includegraphics[width=39mm]{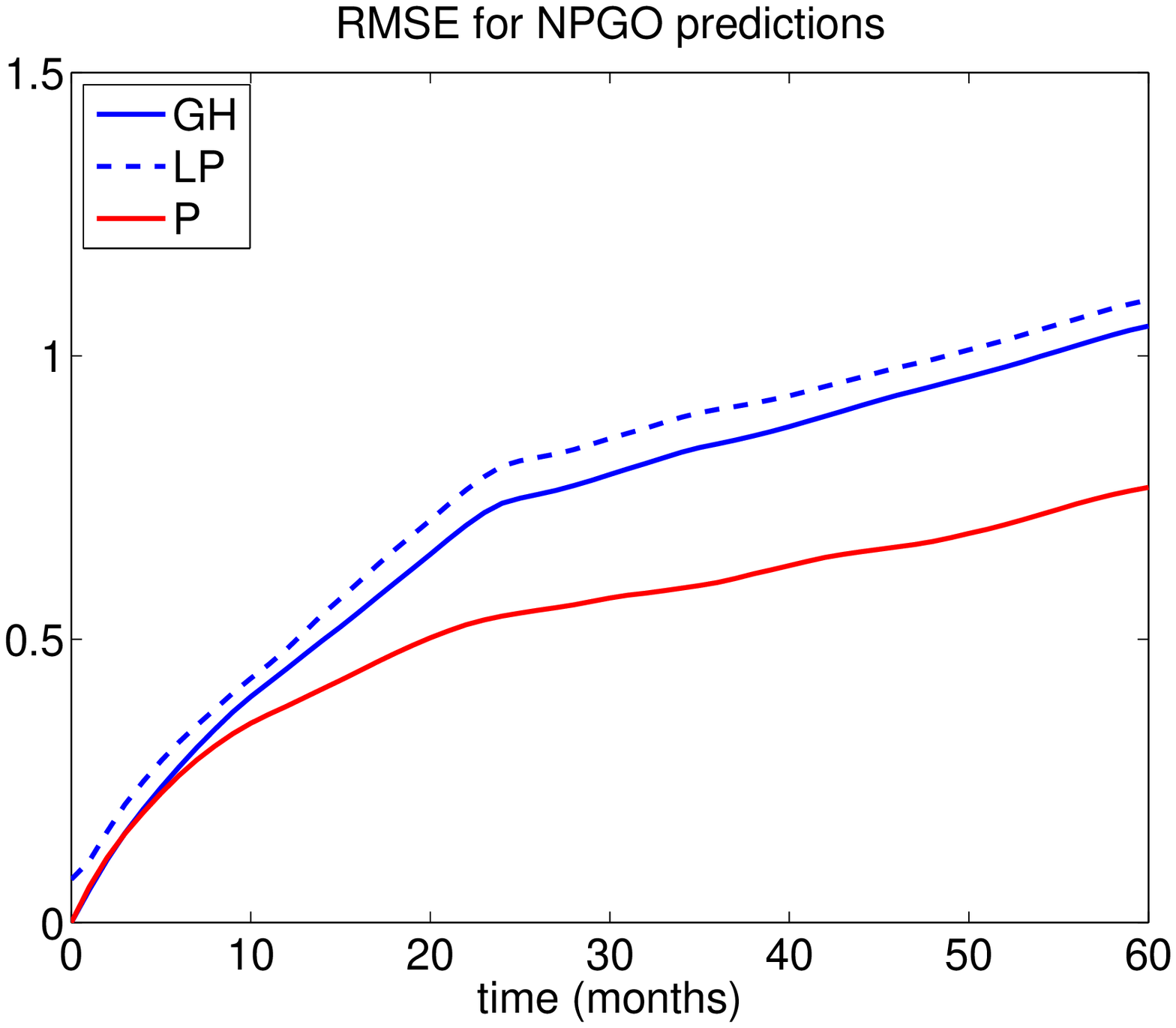}
\includegraphics[width=39mm]{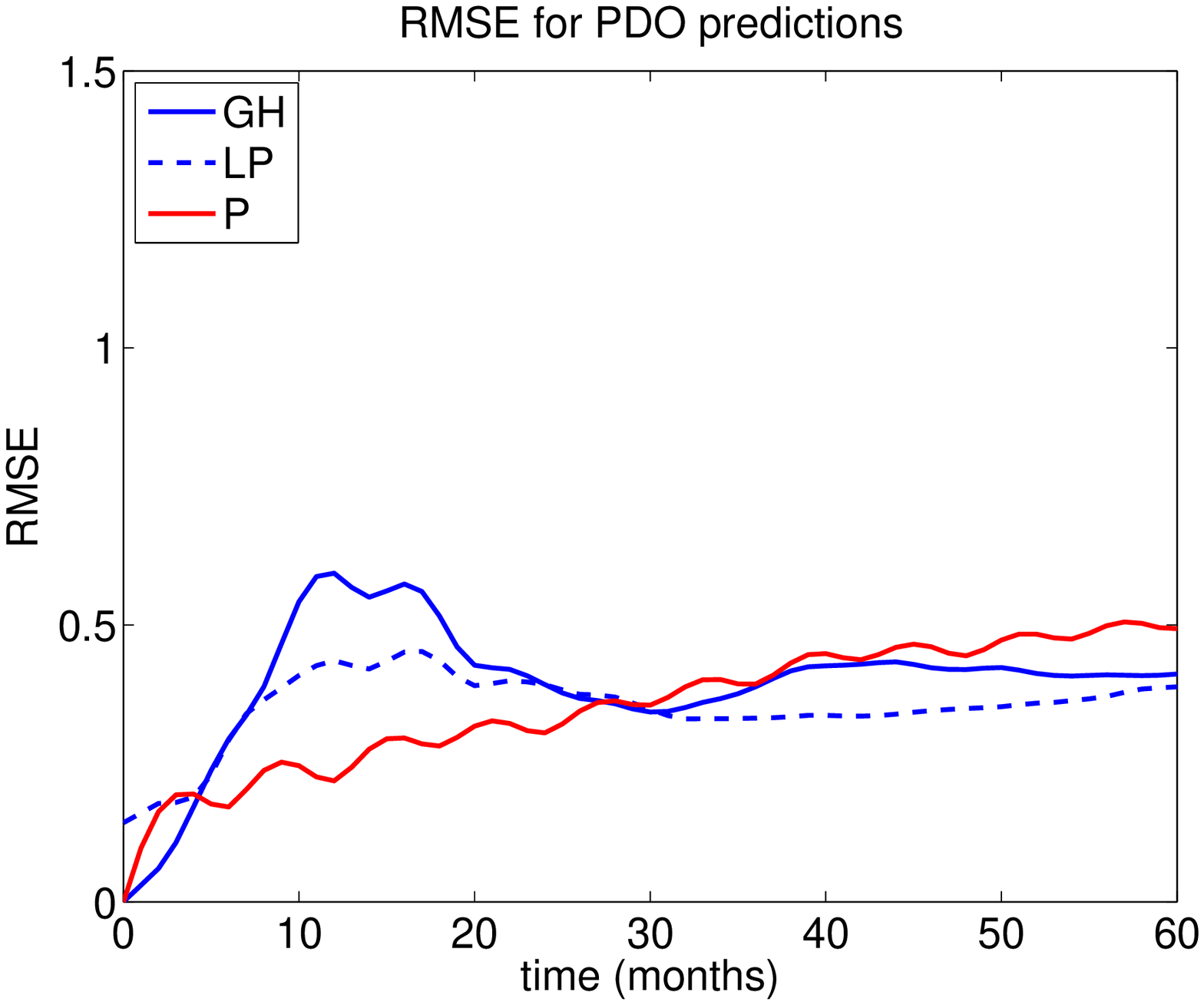}
\includegraphics[width=39mm]{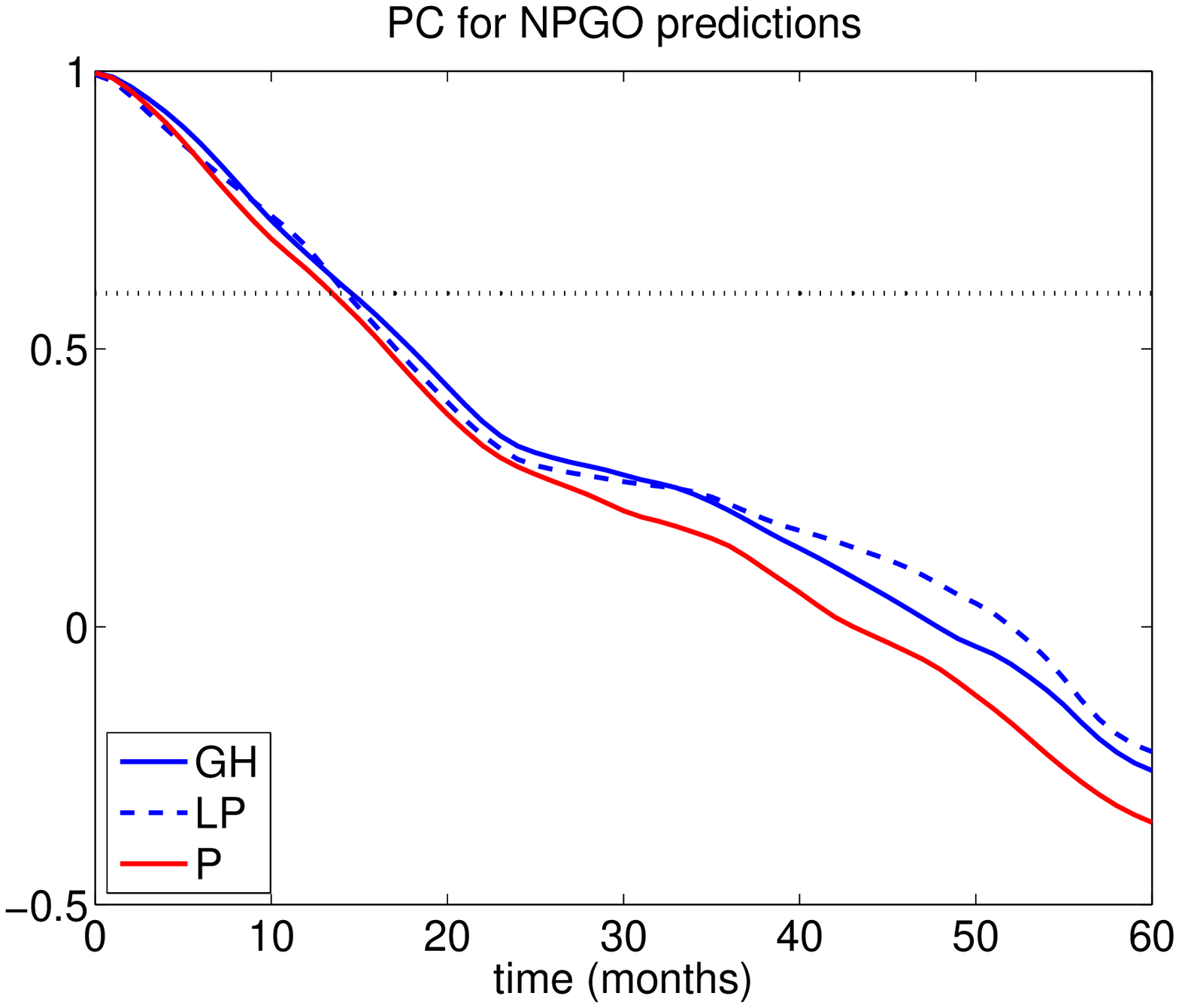}
\includegraphics[width=39mm]{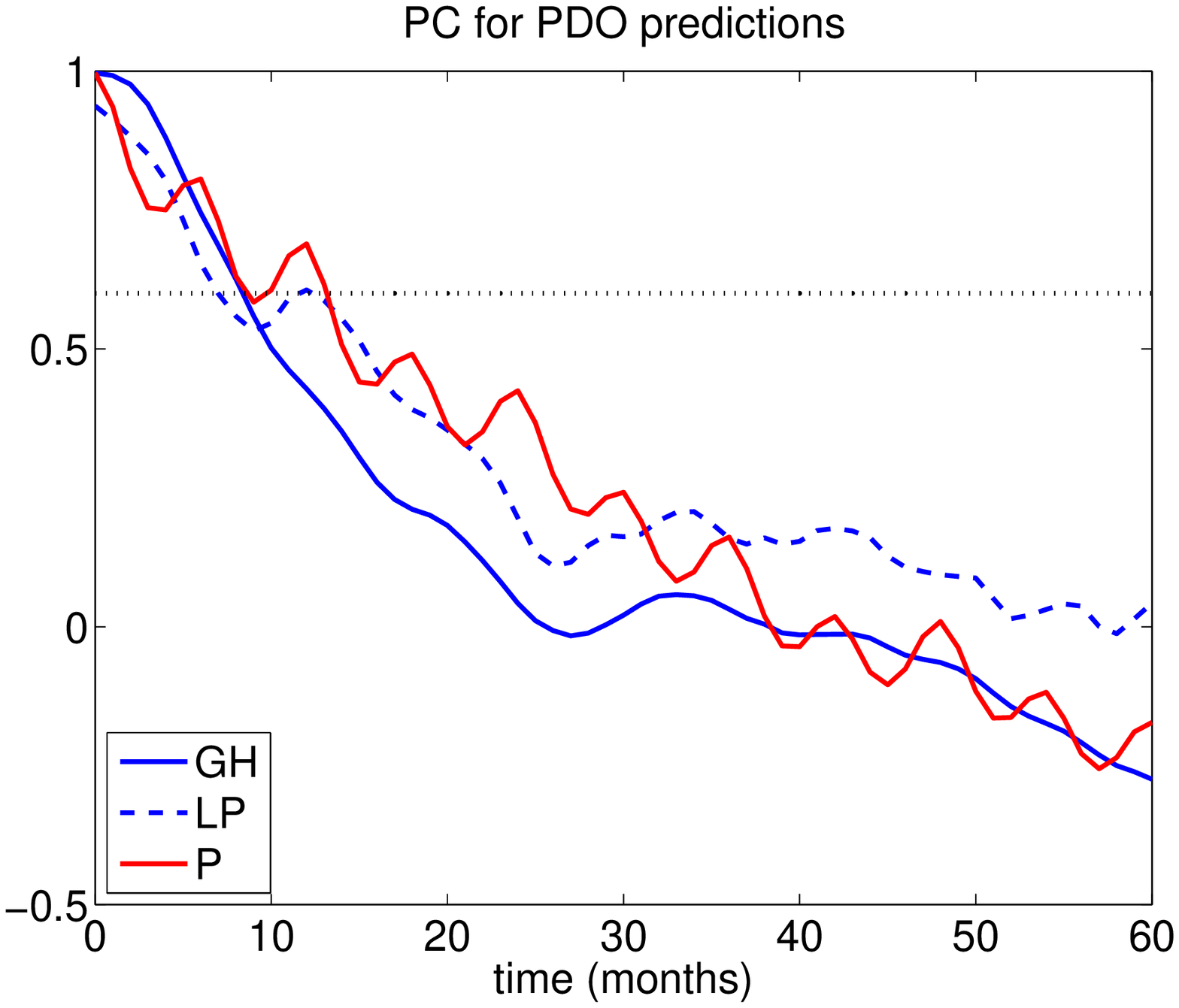}
\caption[KEAF - Model Error w/ Observations - Error Metrics]{Ensemble analog prediction results for the leading two low-frequency modes, with model error. CCSM4 model data is used for the in-sample data, and HadISST observational data is used as the out-of-sample data. Here the method produces little advantage over persistence, given the model error between model and nature, and actually fails to beat persistence for the $\phi^{SI}_{L_2}$ (PDO) mode.}
\label{fig:ccsm4_hadisst_sstice_error}
\end{center}
\end{figure}

 We can further examine the model error scenario by using the actual obervational data set as our nature, and CCSM4 as our model (Figure \ref{fig:ccsm4_hadisst_sstice_error}). Given the short observational record used, far fewer prediction realizations are generated, adding to the noisiness of the error metrics. The loss of predictability is apparent, especially in the $\phi_{L_2}^{SI}$ mode, where the kernel ensemble analog forecasts fail to beat persistence, and drop below 0.6 in PC by 10 months, half as long as the perfect model case.
		
		\subsubsection{Comparison with Autoregressive Models}
		\par
\begin{figure}[ht]
\begin{center}
\includegraphics[width=39mm]{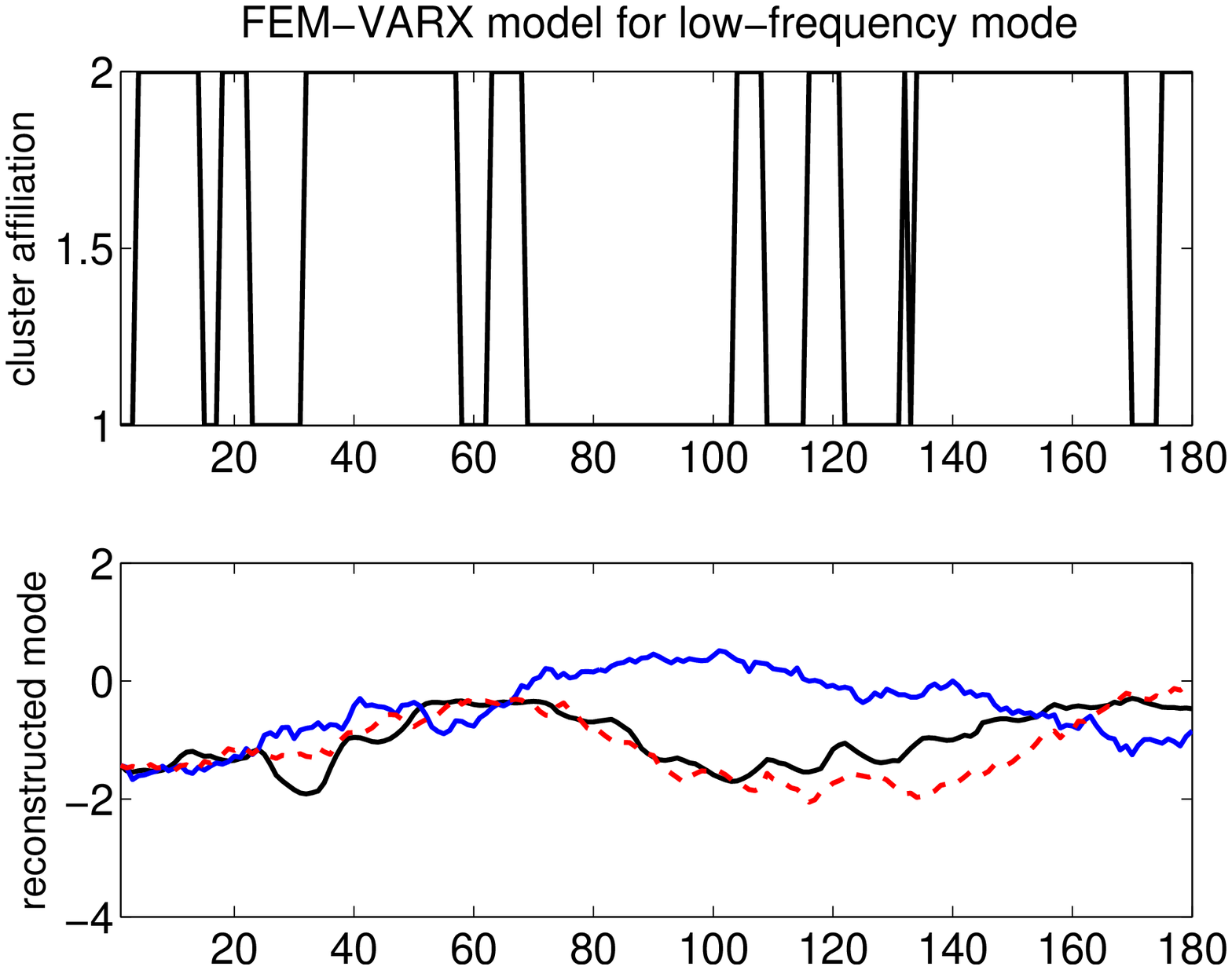}
\includegraphics[width=39mm]{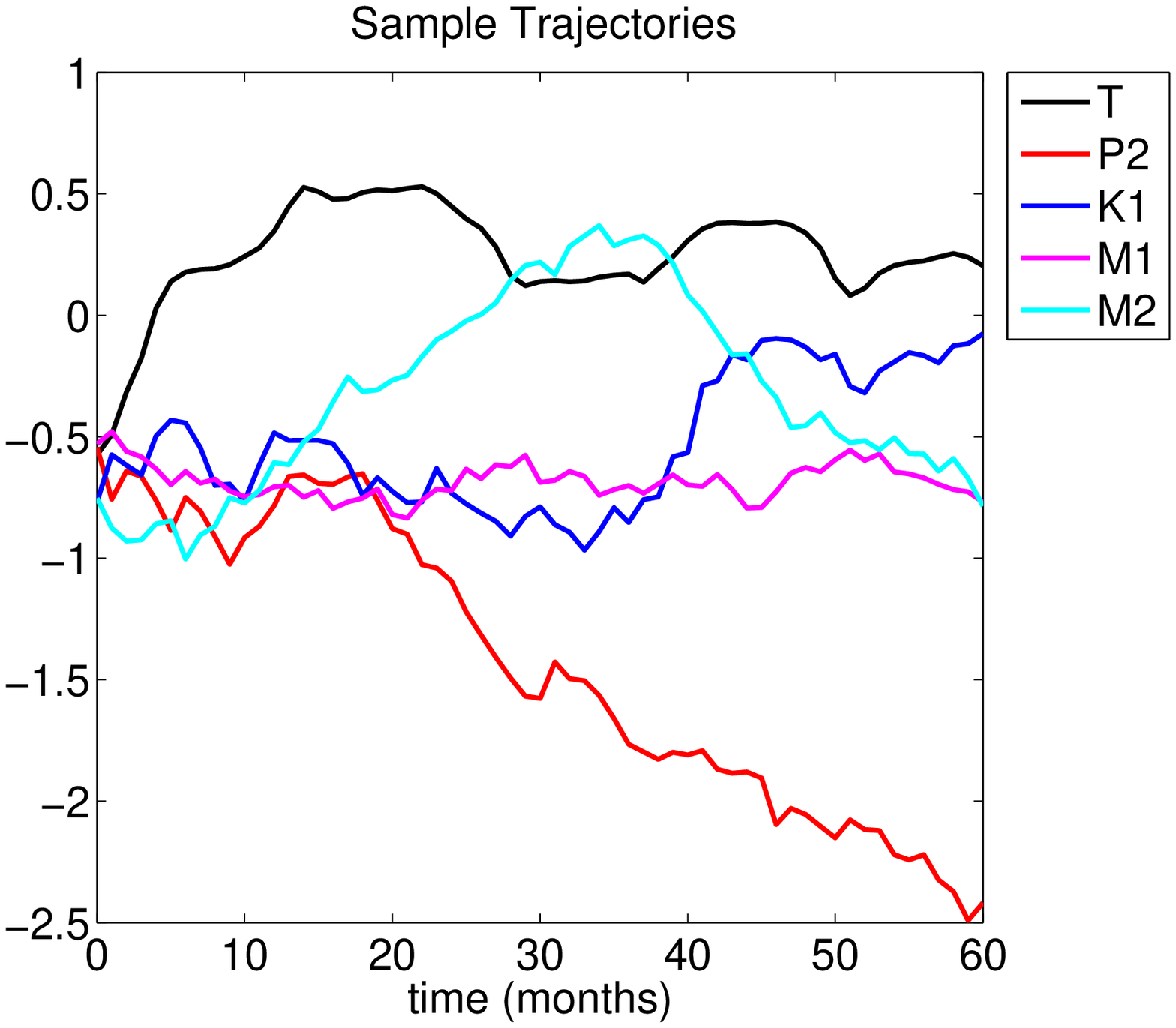}
\includegraphics[width=39mm]{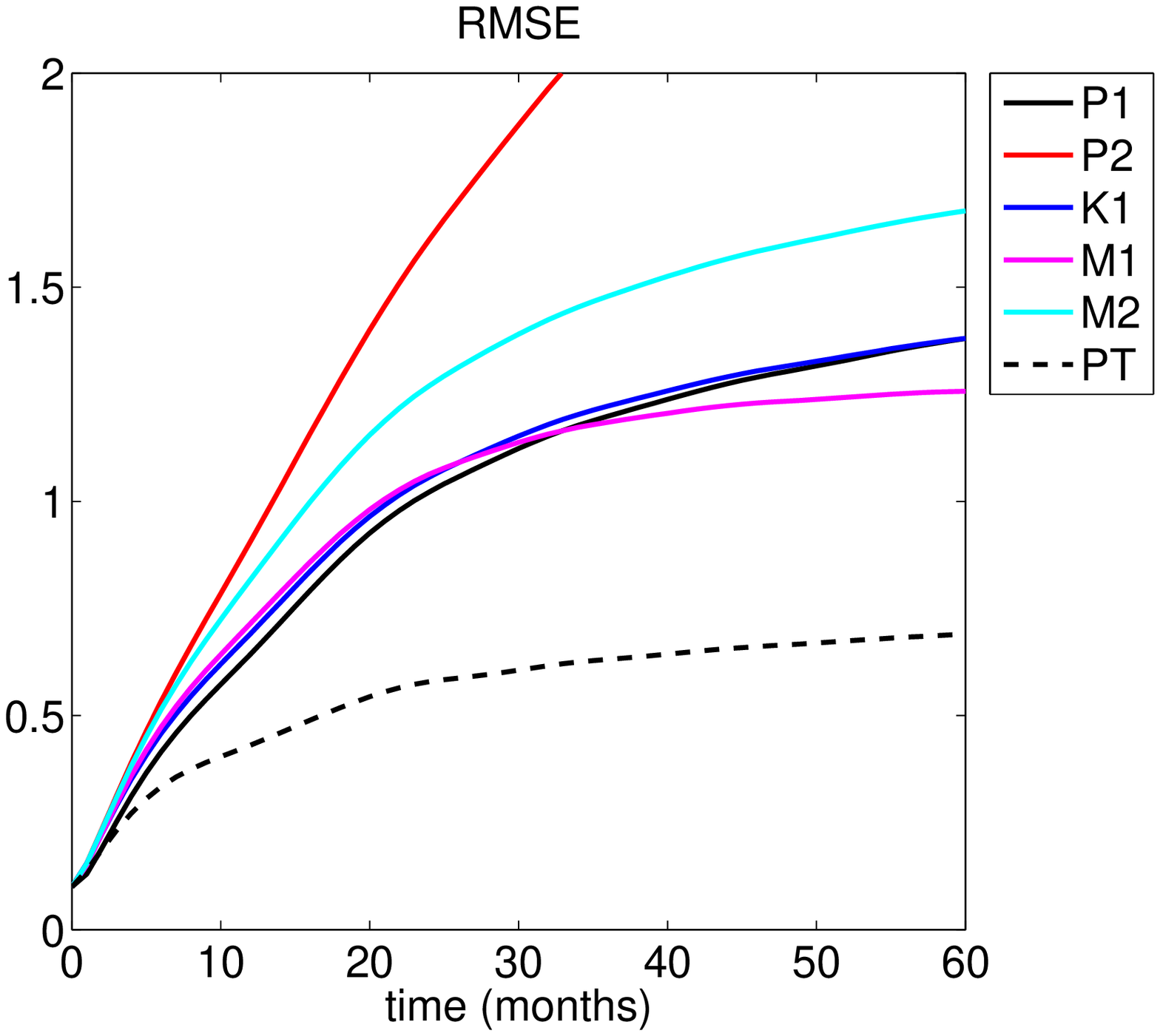}
\includegraphics[width=39mm]{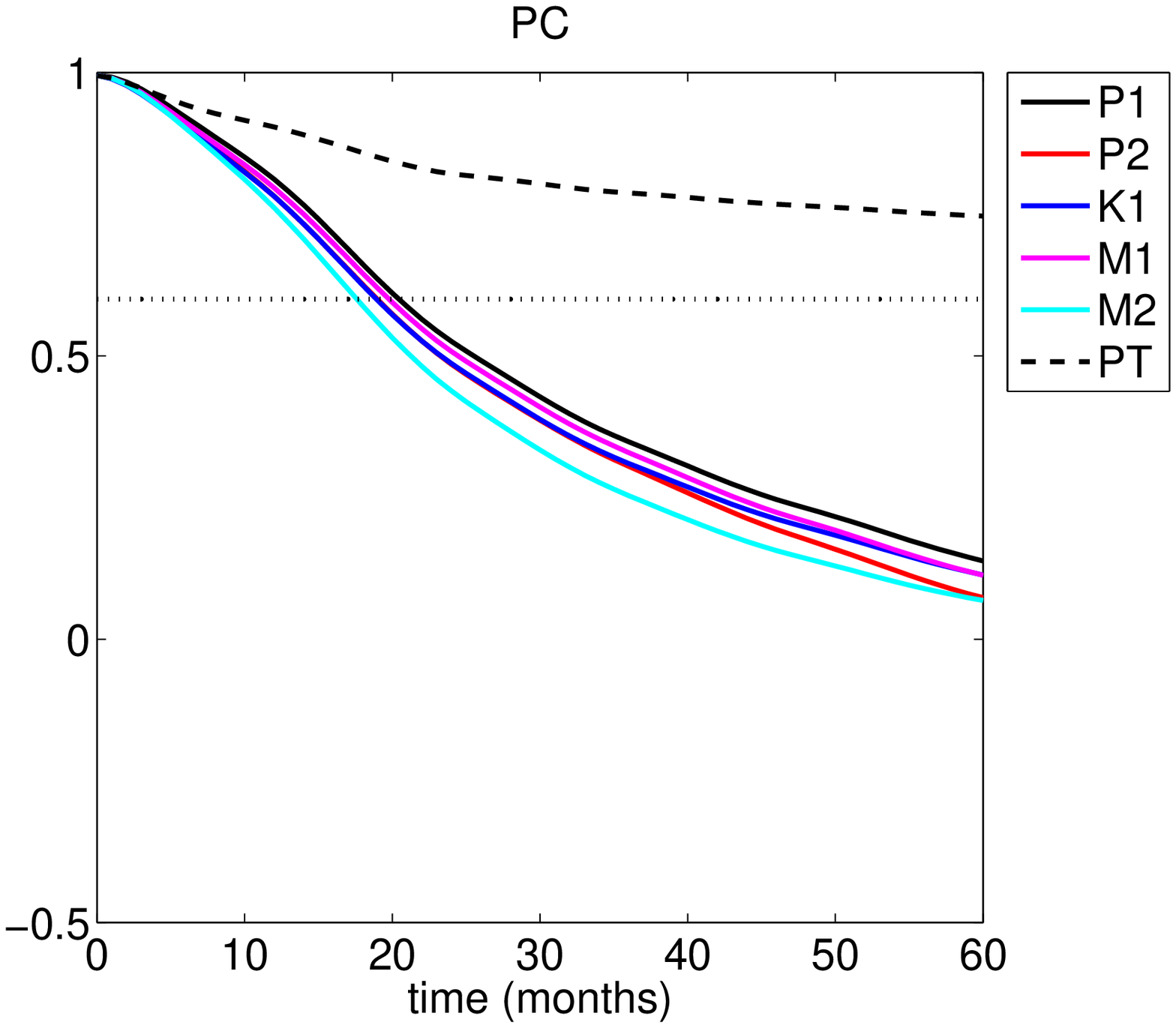}
\caption[AR - CCSM4]{Top left: Snapshot of the true CCSM4 $\phi_{L_1}^{SI}$ (NPGO) trajectory (black) with reconstructed stationary (blue) and non-stationary $K=2$ (red) FEM-VARX model trajectories, along with corresponding model affiliation function $\Gamma(t)$ for non-stationary case. Top right: Sample trajectories for various prediction methods: P2 = stationary, using FEM-VARX model coefficients from initial cluster; K1 = stationary autoregressive; M1, M2 = FEM-VARX with predictions as described in Section \ref{subsubsec:pred_pi}, where M1 is deterministic evolution of the cluster affiliation $\pi(t)$, and M2 uses realizations of $\pi$ generated from the estimated probability transition matrix $T$. Bottom panels: RMSE and PC as a function of lead time for various prediction methods, including P1 = persistence as a benchmark. The dashed black line is for potential predictive skill of non-stationary FEM-VARX, where predictions were ran over the training period using the known optimal model affiliation function $\Gamma(t)$.}
\label{fig:AR_ccsm4_sstice}
\end{center}
\end{figure}

\begin{figure}[ht]
\begin{center}
\includegraphics[width=39mm]{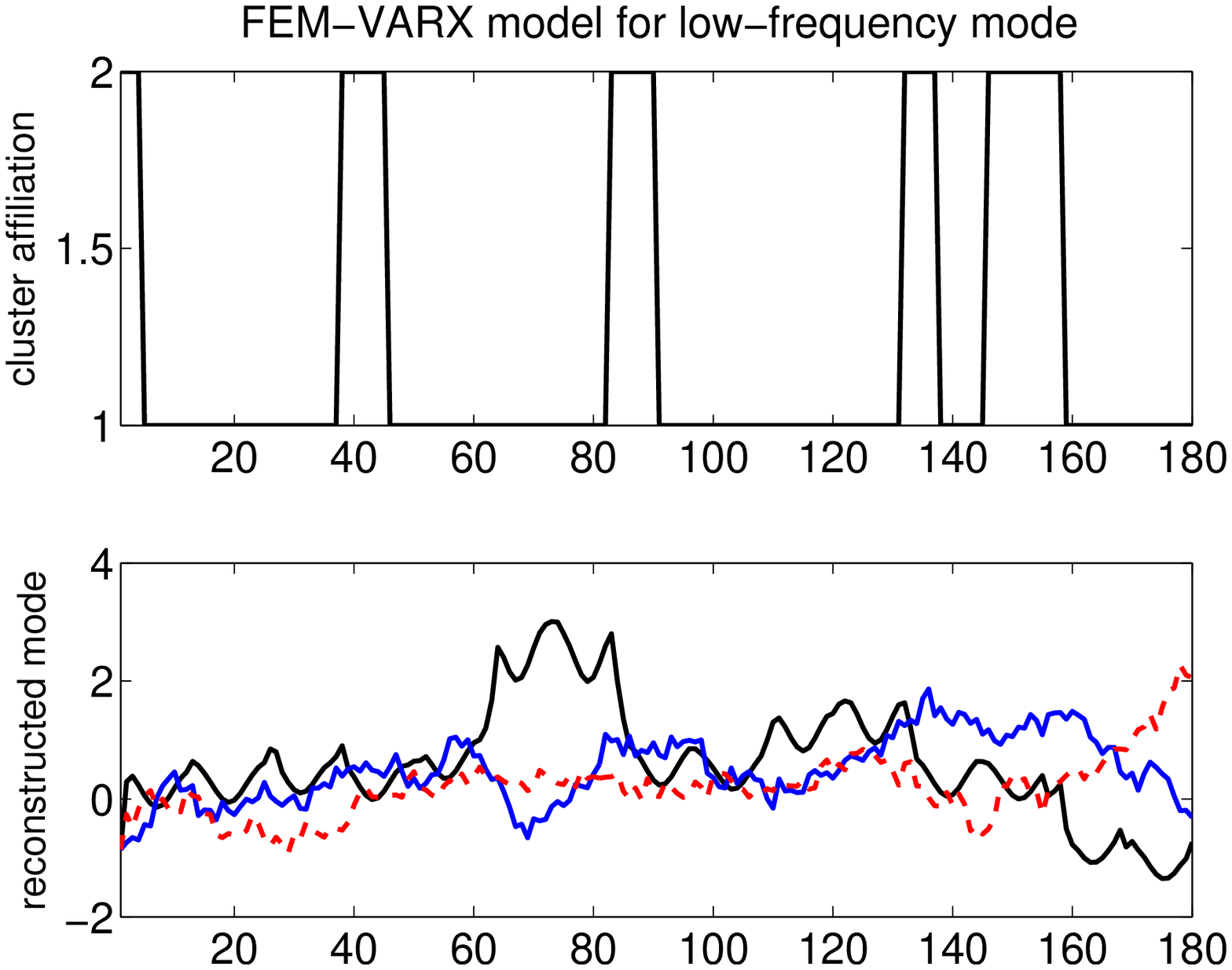}
\includegraphics[width=39mm]{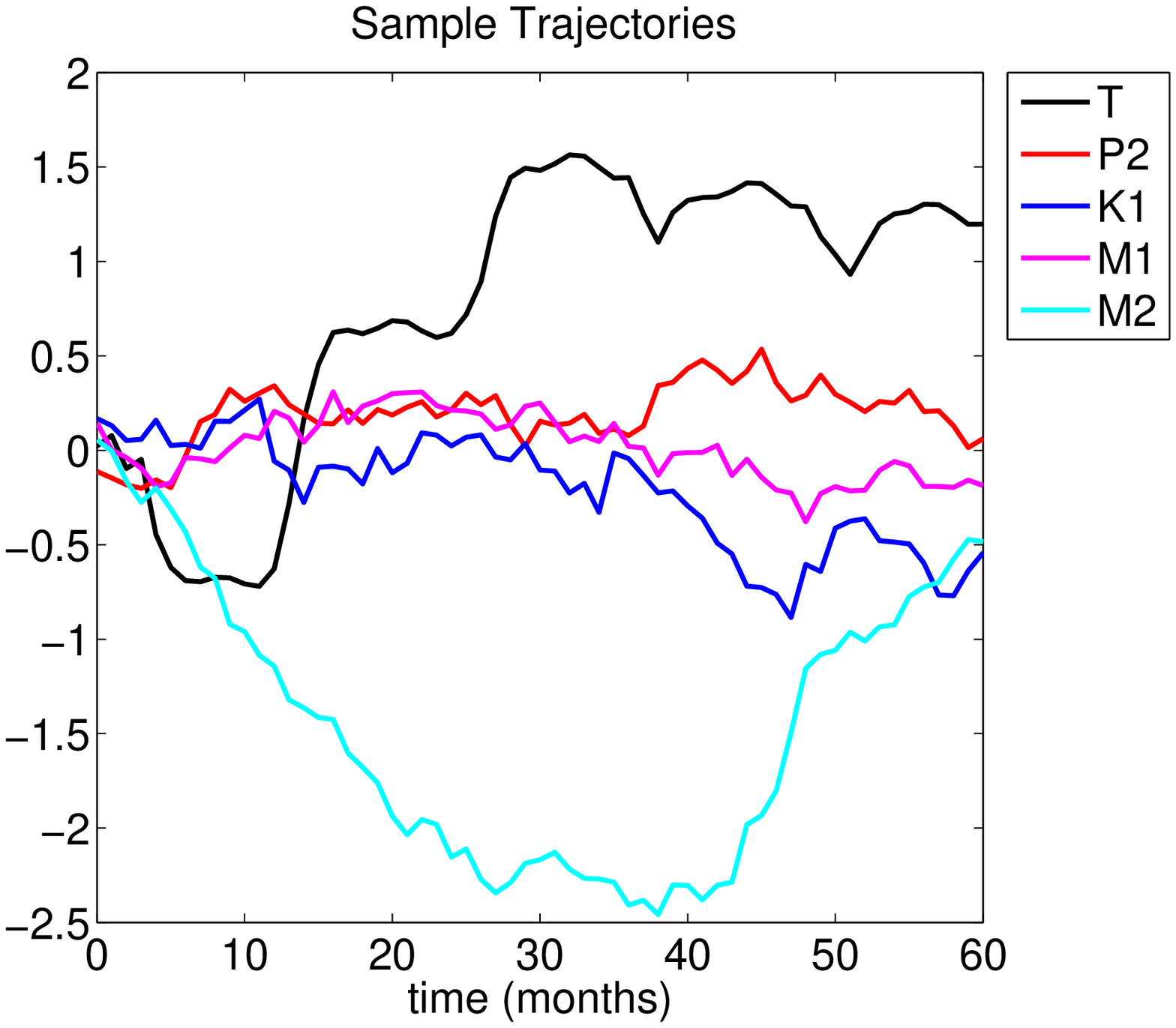}
\includegraphics[width=39mm]{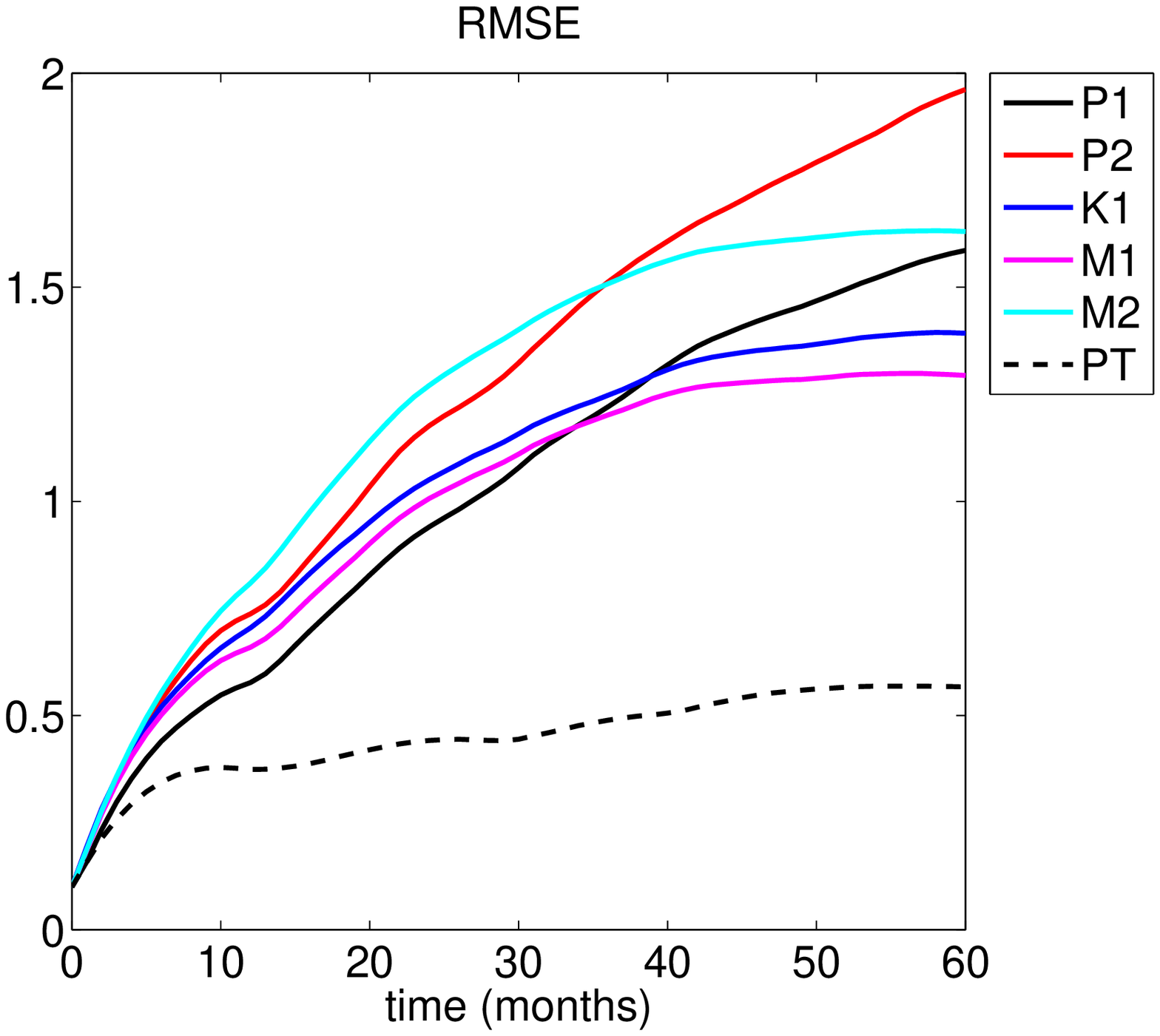}
\includegraphics[width=39mm]{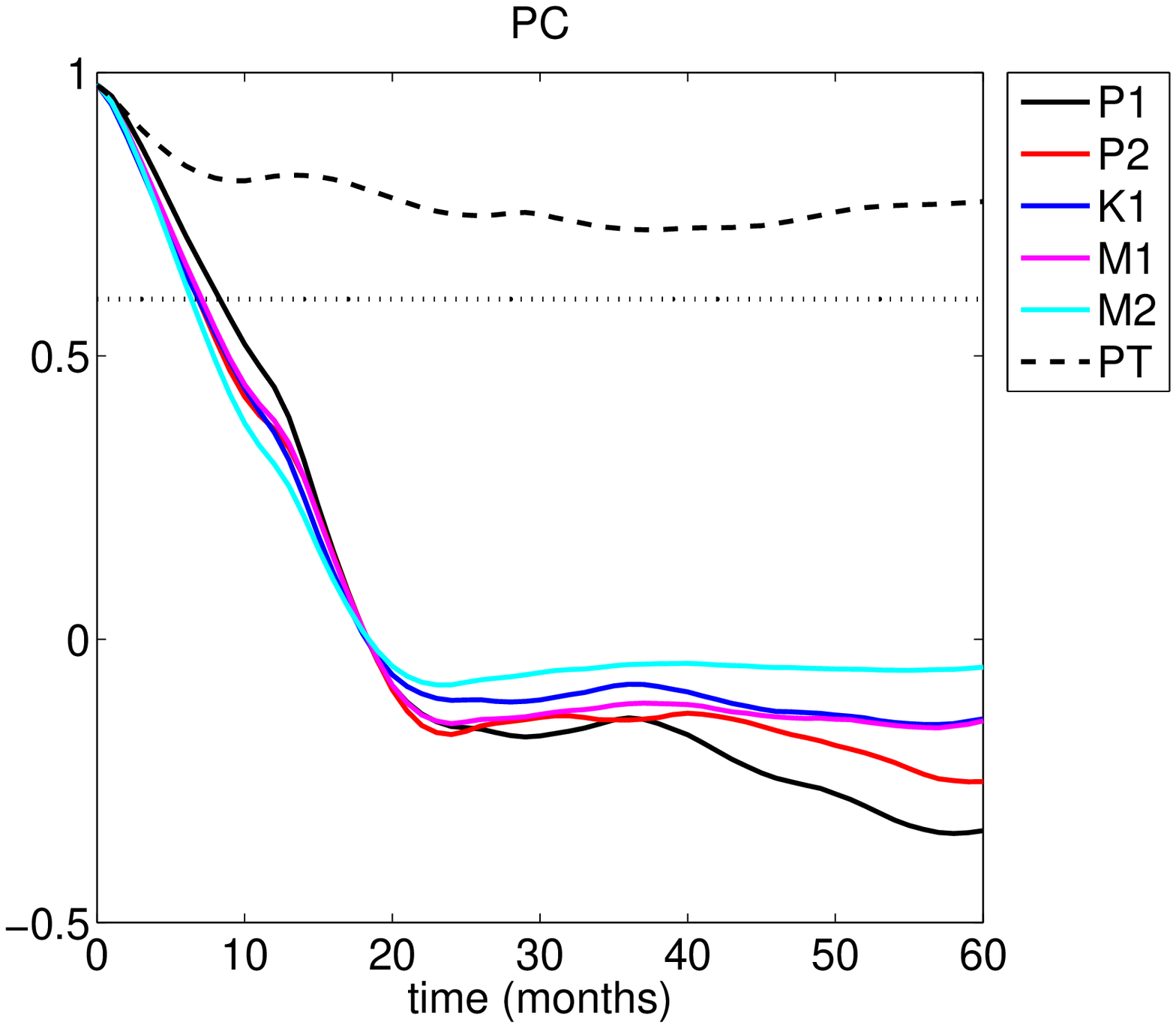}
\caption[AR - HadISST]{Top left: Snapshot of the true HadISST $\phi_{L_1}^{SI}$ (NPGO) trajectory (black) with stationary (blue) and non-stationary $K=2$ (red) FEM-VARX model trajectories, along with the corresponding model affiliation function $\Gamma(t)$ for non-stationary case. Top right: Sample trajectories for various prediction methods--see Figure \ref{fig:AR_ccsm4_sstice} for details of methods. Bottom panels: RMSE and PC as a function of lead time for various prediction methods. The dashed black line is for potential predictive skill of non-stationary FEM-VARX, where predictions were ran over the training period using the known optimal model affiliation function $\Gamma(t)$.}
\label{fig:AR_hadisst_sstice}
\end{center}
\end{figure}

 We compare the ensemble analog predictions to standard stationary autoregressive models, as well as non-stationary models using the FEM-VARX framework of \citet{horenko2010identification} discussed in Section \ref{subsec:ARmodel}, for the low-frequency modes $\phi^{SI}_{L_1}$, $\phi^{SI}_{L_2}$ generated from CCSM4 model (Figure \ref{fig:AR_ccsm4_sstice}) and the HadISST observation (\ref{fig:AR_hadisst_sstice}) training data. For both sets of low-frequency modes, $K=2$ clusters was judged to be optimal by the AIC as mentioned in Section \ref{subsec:ARmodel}--see \citet{horenko2010identification,metzner2012analysis} for more details. The coefficients for each cluster are nearly the same (autoregressive coefficent close to 1, similar noise coefficients), apart from the constant forcing coefficient $\mu_i$ of almost equal magnitude and opposite sign, suggesting two distinct regime behaviors. In the stationary case, the external forcing coefficient $\mu$ is very close to 0.
		\par In the top left panel of Figures \ref{fig:AR_ccsm4_sstice} and \ref{fig:AR_hadisst_sstice}, we display snapshots of the leading low-frequency mode $\phi_{L_1}^{SI}$ (NPGO) trajectory reconstruction during the training period, for the stationary (blue, solid), and non-stationary (red, dashed) models, along with the cluster switching function associated with the non-stationary model. In both the CCSM4 model and HadISST data sets, the non-stationary model snapshot is a better representation of the truth (black, solid), and the benefit over the stationary model is more clearly seen in the CCSM4 model data, Figure \ref{fig:AR_ccsm4_sstice}, which has the benefit of a 400 year training period, as opposed to the shorter 16 year training period with the observational data set.
		\par In the prediction setting, however, the non-stationary models, which are reliant on an advancement of the unknown cluster affiliation function $\Gamma(t)$ beyond the training period, as discussed in Section \ref{subsec:ARmodel}, fail to outperform their stationary counterparts in the RMSE and PC metrics (bottom panels of Figures \ref{fig:AR_ccsm4_sstice} and \ref{fig:AR_hadisst_sstice}). In fact, none of the proposed regression prediction models are able to outperform the simple persistence forecast in these experiments. As a measure of \emph{potential} predition skill for the non-stationary models, whereby we mean that if perfect knowledge of the underlying optimal cluster switching function $\Gamma(t)$ could be known over the test period, we have run the experiment of replacing the test data period with the training data set, and find exceedingly strong predictive performance, with PC between 0.7 and 0.8 for all time lags tested, up to 60 months. Similar qualitative results for the second leading low-frequency mode $\phi_{L_2}^{SI}$ (PDO) for each data set were found (not shown). This suggests that the Markov hypothesis, the basis for the predictions of $\Gamma(t)$, is not accurate, and other methods incorporating more memory are needed.
		
	\subsection{Sea ice anomalies}
	\label{subsec:siea}
	\par
\begin{figure}[ht]
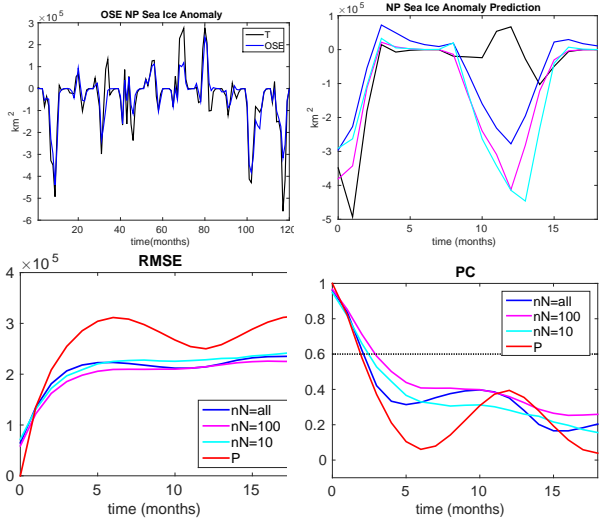

\begin{center}
\includegraphics[width=39mm]{Fig9a}
\includegraphics[width=39mm]{Fig9b}
\includegraphics[width=39mm]{Fig9c}
\includegraphics[width=39mm]{Fig9d}
\caption[North Pacific Sea Ice Cover Anomalies]{Top left: True sea ice cover anomalies plotted with the corresponding Laplacian pyramid out-of-sample extention function. The other panels show prediction results using Laplacian pyramids for total North Pacific sea ice cover anomalies in CCSM4 data. The number of nearest neighbors (nN) used to form the ensemble was varied, and we find the best performance when the ensemble is restricted to the nearest 100 neighbors, corresponding to about 2\% of the total sample size.}
\label{fig:ccsm4_ccsm4_NPICA}
\end{center}
\end{figure}

 The targeted observables hereto considered for prediction have been data driven, and as such influenced by the data analysis algorithm. Hence there is no objective ground truth available when predicting these modes beyond the training period on which the data analysis was performed, and while in this case the NLSA algorithm was used, other data analysis methods such as EOFs would suffer the same drawback. We wish to test our prediction method on an observable that is objective, in the sense that it can be computed independently of the data analysis algorithm, for which we turn to integrated sea ice extent anomalies, as defined in Equation \eqref{eq:siea}. We can clearly compute the time series of sea ice anomalies from the out-of-sample set directly (relative to the training set climatology), which will be our ground truth, and use the Laplacian pyramid approach to generate our out-of-sample extension predictions. This observable does not have a tight expansion in the eigenfunction basis, so the geometric harmonics method of extension will be ill-conditioned, and thus not considered. We note that in this approach, there are reconstruction errors at time lag $\tau=0$, so at very short time scales we cannot outperform persistence. We consider a range of truncation levels for the number of ensemble members used, which are nearest neighbors to the out-of-sample data point, as determined by the kernel function. Using all available neighbors will likely overly smooth and average out features in forward trajectories, while using too few neighbors will place to much weight on particular trajectories. Indeed we find good performance using 100 (out of total possible 4791). The top left panel of Figure \ref{fig:ccsm4_ccsm4_NPICA} shows a snapshot of the true sea ice extent anomalies, respectively, together with a reconstructed out-of-sample extension using the Laplacian pyramid. To be clear this is not a prediction trajectory, but rather each point in the out-of-sample extension is calculated using Equation \eqref{eq:LP}; that is, each point is a time-lead $\tau=0$ reconstruction.
	
	 \par The top right panel shows sample snapshots of prediction trajectories, restricting the ensemble size to the nearest 10, 100, and then all nearest neighbors. Notice in particular predictions match the truth when the anomalies are close to 0, but then may subsequently progress in the opposite sign as the truth. As our prediction metrics are averaged over initial conditions spanning all months, the difficulty the predictions have in projecting from a state of near 0 anomaly significantly hampers the ability for long-range predictability of this observable.
	 \par In the bottom two panels of Figure \ref{fig:ccsm4_ccsm4_NPICA} we have the averaged error metrics, and see year to year correlations manifesting as a dip/bump in RMSE and PC in the persistence forecast that occurs after 12 months. After the first month lag time, the kernel ensemble analog forecasts overcome the reconstruction error and beat persistence in both RMSE and PC, and give about a 2 month increase in prediction skill (as measured by when the PC drops below 0.6) over persistence. We see the best performance restricting the ensemble size to 100 nearest neighbors (about 2\% of the total sample size) in both the RMSE and PC metrics, though this is marginal before the error metrics drop below the 0.6 threshold.
	 \par
\begin{figure}[ht]
\begin{center}
\includegraphics[width=39mm]{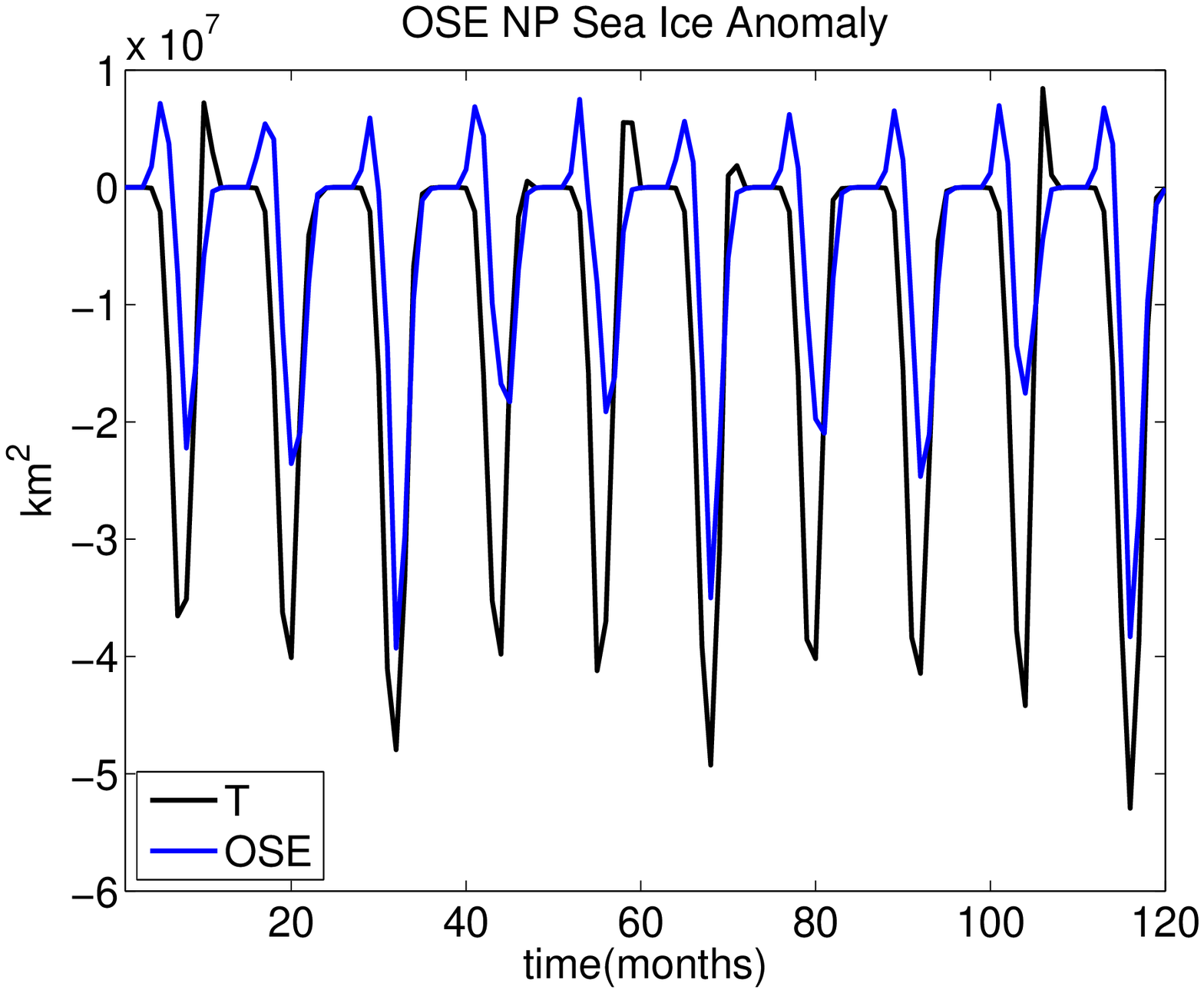}
\includegraphics[width=39mm]{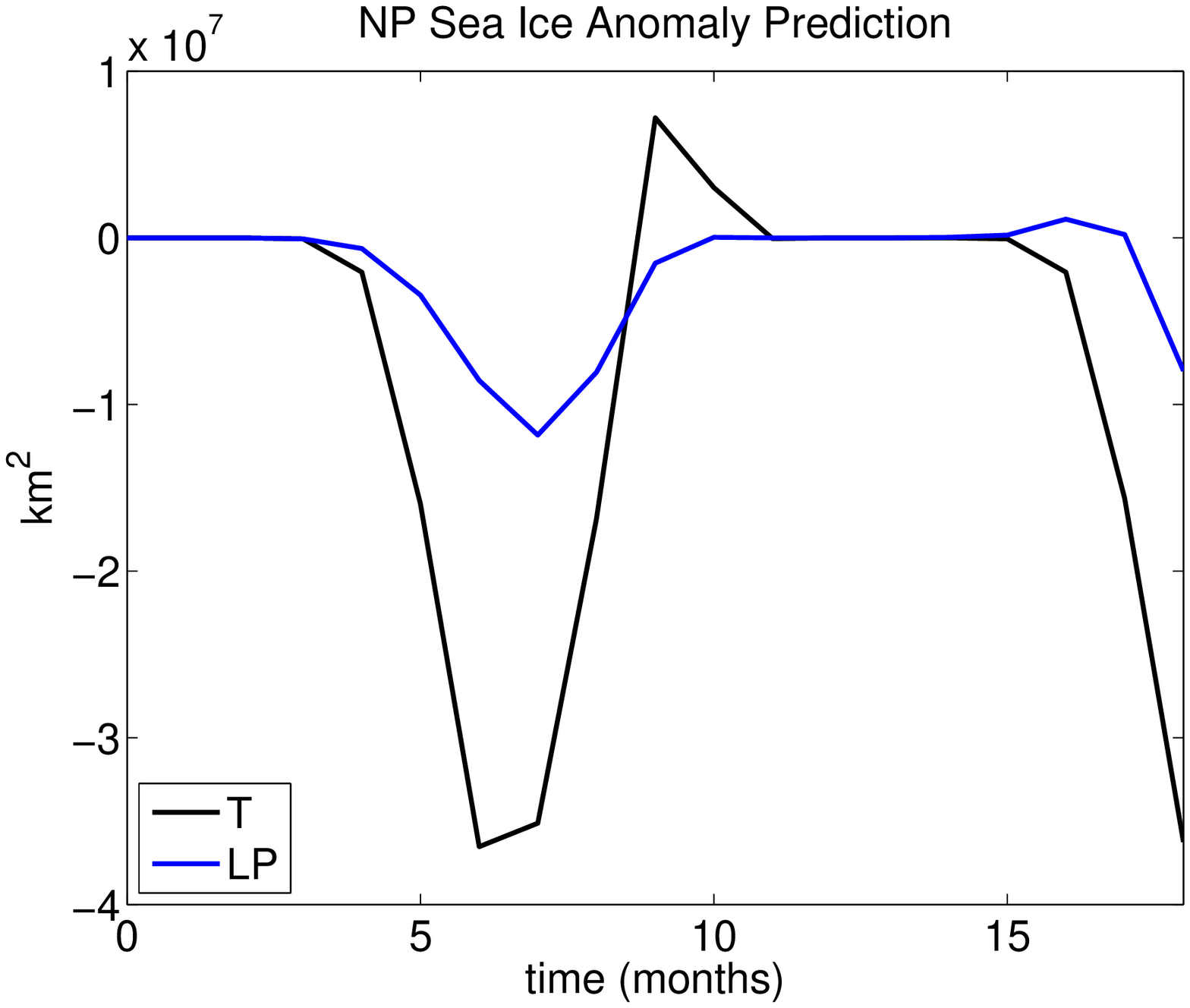}
\includegraphics[width=39mm]{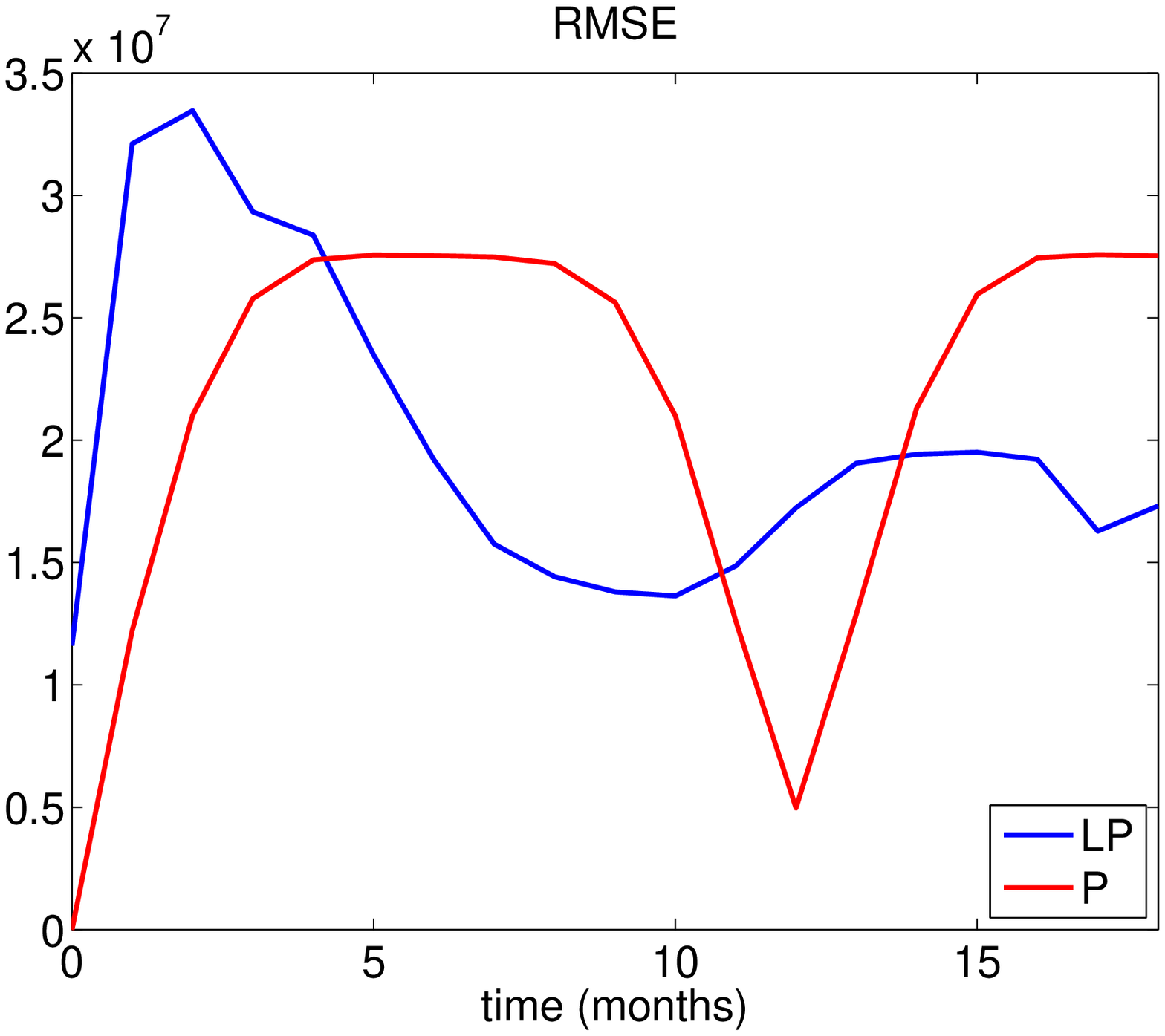}
\includegraphics[width=39mm]{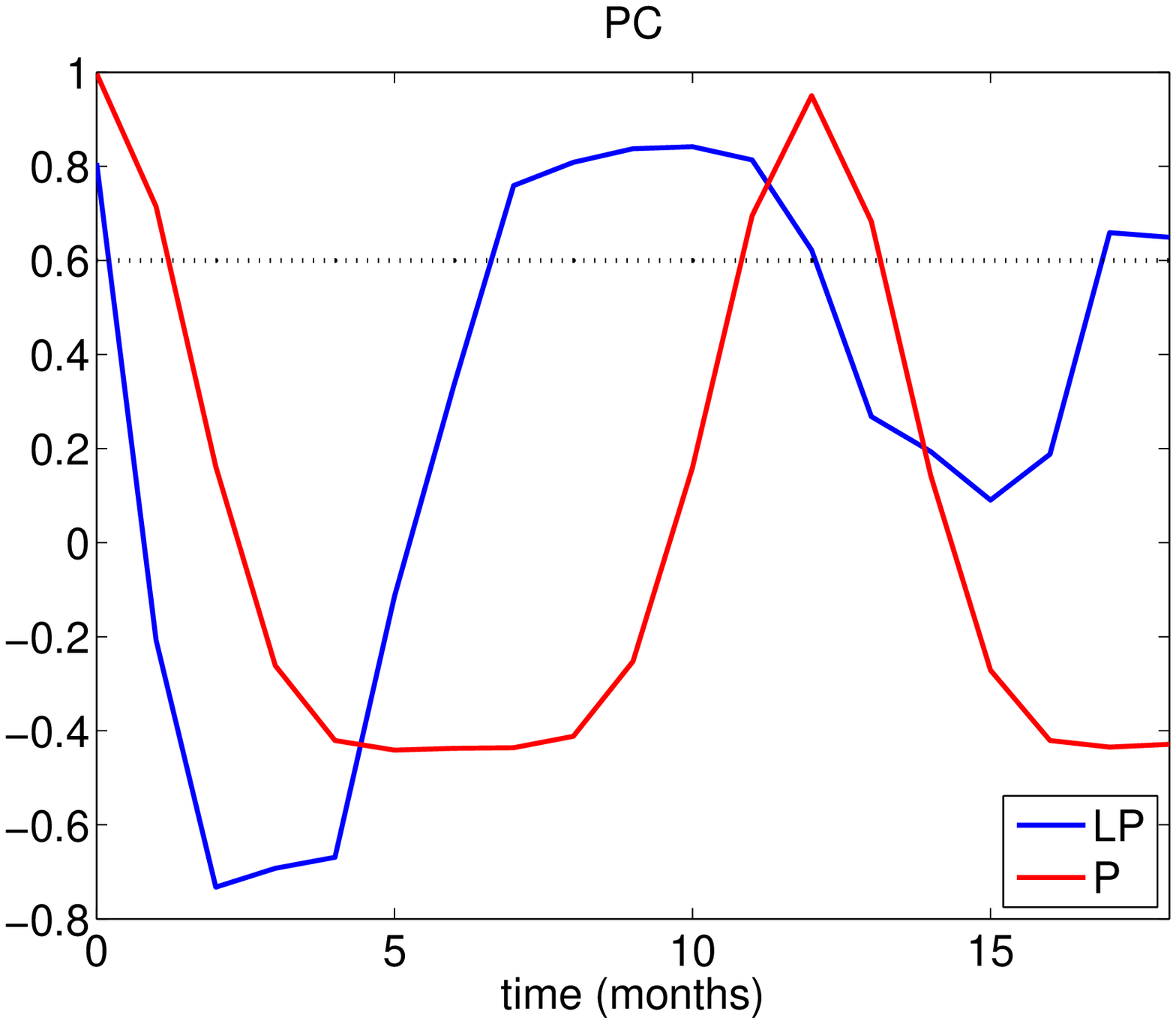}
\caption[North Pacific Sea Ice Volume Anomalies]{Prediction results for total observed (HadISST data) North Pacific sea ice volume anomalies, using CCSM4 as training data. The ice cover function is out-of-sample extended via Laplacian pyramid, using 100 nearest neighbors.}
\label{fig:ccsm4_hadisst_NPICA}
\end{center}
\end{figure}

 Pushing the prediction strategy to an even more difficult problem, in Figure \ref{fig:ccsm4_hadisst_NPICA} we try to predict observational sea ice extent anomalies using CCSM4 model data as training data. In this scenario, without knowledge of the test data climatology, the observation sea ice extent anomalies are defined using the CCSM4 climatology. In the top panel of Figure \ref{fig:ccsm4_hadisst_NPICA} we see the strong bias as a result, where the observational record has less sea ice than the CCSM model climatology, which has been taken from a pre-industrial control run. This strongly hampers the ability to accurately predict observational sea ice extent anomalies using CCSM4 model ensemble analogs, and as a result the only predictive skill we see is from the annual cycle. 
	 
\section{Discussion}
\label{sec:discussion}
	\par We have examined a recently proposed prediction strategy employing a kernel ensemble analog forecasting scheme making use of out-of-sample extension techniques. These nonparametric, data-driven methods make no assumptions on the underlying governing dynamics or statistics. We have used these methods in conjunction with NLSA to extract low-frequency modes of variability from North Pacific SST and SIC data sets, both from models and observations. We find that for these low-frequency modes, the analog forecasting performs at least as well, and in many cases better than, the simple constant persistence forecast. Predictive skill, as measured by PC exceeding 0.6, can be increased by up to 3 to 6 months for low-frequency modes of variability in the North Pacific. This is a strong advantage over traditional parametric regression models, which were shown to fail to beat persistence.
	\par The kernel ensemble analog forecasting methods outlined included two variations on the underlying out-of-sample extension scheme, each with its strengths and weaknesses. The geometric harmonics method, based on the Nystr{\"o}m method, worked well for observables that are band-limited in the eigenfunction basis, in particular the eigenfunctions themselves. However for observables not easily expressed in such a basis, the Laplacian pyramid provides an alternative method based on a multiscale decomposition of the original observable.
	\par While the low-frequency eigenfunctions from NLSA were a natural preferred class of observables to target for prediction, we also studied the case of objective observables uninfluenced by the data analysis algorithm. Motivated by the strong reemergence phenomena, we considered sea ice extent anomalies as our target for prediction in the North Pacific. Using a shorter embedding window due to faster (seasonal) time scale dynamics, we obtain approximately a two month increase in predictive skill over the persistence forecast. It is also evident that when considering regional sea ice extent anomalies winds play a large role in moving ice into and out of the domain of interest, and as such additional consideration of the atmospheric component in the system could be included in the multivariate kernel function, despite having weaker low-frequency variability.
	\par An important consideration is that our prediction metrics are averaged over initial conditions ranging over all possible initial states of the system. As we saw clearly in the case of North Pacific sea ice volume anomalies, these prediction strategies can have difficulty with projecting from an initial state of quiessence, and can easily predict to the wrong sign of an active state, greatly hampering predictive skill. On the other hand we would expect predictive skill to be stronger for those initial states that begin in a strongly active state, or said differently, clearly in one climate regime, as oppose to in transition between the two. Future work will further explore conditional forecasting, where we either condition forecasts on the initial month, or the target month. Also extending this analysis to the North Atlantic, another region of strong low-frequency variability, is a natural progression of this work.




\begin{acknowledgements}
The research of Andrew Majda and Dimitrios Giannakis is partially supported by ONR MURI grant 25-74200-F7112. Darin Comeau is supported as a postdoctoral fellow through this grant. The research of Dimitrios Giannakis is also partially supported by ONR DRI grant N00014-14-0150.

\end{acknowledgements}



\end{document}